\documentclass[10pt]{article}

\usepackage{amsmath}
\usepackage{amsthm}
\usepackage{amssymb}
\usepackage{amscd}
\usepackage{amsfonts}  
\usepackage{makeidx}
\usepackage[all]{xy} \SelectTips{cm}{}
\usepackage{makeidx}     
\usepackage{graphicx}    
\usepackage{multicol}    

\DeclareMathAlphabet\EuR{U}{eur}{m}{n}
\SetMathAlphabet\EuR{bold}{U}{eur}{b}{n}

\makeindex             

\theoremstyle{plain}
\newtheorem{theorem}{Theorem}[section]
\newtheorem{lemma}[theorem]{Lemma}
\newtheorem{proposition}[theorem]{Proposition}
\newtheorem{corollary}[theorem]{Corollary}

\theoremstyle{definition}
\newtheorem{definition}[theorem]{Definition}
\newtheorem{example}[theorem]{Example}
\newtheorem{remark}[theorem]{Remark}
\newtheorem{notation}[theorem]{Notation}

{\catcode`@=11\global\let\c@equation=\c@theorem}

\newcommand{\comsquare}[8]                   
{\begin{CD}
#1 @>#2>> #3\\
@V{#4}VV @V{#5}VV\\
#6 @>#7>> #8
\end{CD}
}

\newcommand{\xycomsquare}[8]                   
{\xymatrix
{#1 \ar[r]^{#2} \ar[d]^{#4} &
#3 \ar[d]^{#5}  \\
#6\ar[r]^{#7} &
#8
}
}

\newcommand{\calall}{\mathcal{ALL}}

\newcommand{\calfin}{\mathcal{FIN}}
\newcommand{\calsfg}{\mathcal{SFG}}
\newcommand{\calsub}{\mathcal{SUB}}
\newcommand{\calvcyc}{\mathcal{VCY}}
\newcommand{\caltr}{\{ \! 1 \! \}}

\newcommand{\calc}{{\cal C}}

\newcommand{\calf}{{\cal F}}
\newcommand{\calg}{{\cal G}}
\newcommand{\calh}{{\cal H}}

\newcommand{\calm}{{\cal M}}

\newcommand{\IN}{{\mathbb N}}

\newcommand{\IQ}{{\mathbb Q}}

\newcommand{\IZ}{{\mathbb Z}}

\newcommand{\bfK}{{\mathbf K}}

\newcommand{\bfL}{{\mathbf L}}

\newcommand{\curs}{\EuR}

\newcommand{\Or}{\curs{Or}}

\newcommand{\SPACES}{\curs{SPACES}}

\newcommand{\aut}{\operatorname{aut}}

\newcommand{\cd}{\operatorname{cd}}
\newcommand{\cent}{\operatorname{cent}}

\newcommand{\coker}{\operatorname{coker}}
\newcommand{\colim}{\operatorname{colim}}

\newcommand{\cone}{\operatorname{cone}}

\newcommand{\cyl}{\operatorname{cyl}}

\newcommand{\hdim}{\operatorname{hdim}}
\newcommand{\Hei}{\operatorname{Hei}}

\newcommand{\id}{\operatorname{id}}

\newcommand{\im}{\operatorname{im}}

\newcommand{\ind}{\operatorname{ind}}

\newcommand{\mor}{\operatorname{mor}}

\newcommand{\ob}{\operatorname{ob}}
\newcommand{\pr}{\operatorname{pr}}

\newcommand{\res}{\operatorname{res}}

\newcommand{\tors}{\operatorname{tors}}

\newcommand{\vcd}{\operatorname{vcd}}

\newcommand{\pt}{\{\bullet\}}

\newcommand{\EGF}[2]{E_{#2}(#1)}                   
\newcommand{\eub}[1]{\underline{E}#1}              
\newcommand{\edub}[1]{\underline{\underline{E}}#1} 
\newcommand{\OrGF}[2]{\Or_{#2}(#1)}                

\newcommand{\higherlim}[3]{{\setbox1=\hbox{\rm lim}
        \setbox2=\hbox to \wd1{\leftarrowfill} \ht2=0pt \dp2=-1pt
        \mathop{\vtop{\baselineskip=5pt\box1\box2}}
        _{#1}}^{#2}#3}

\newcommand{\version}[1]                       
{\begin{center} last edited on #1\\
last compiled on \today
\end{center}
}

\newcounter{commentcounter}

\begin{document}

\typeout{----------------------------  evcyc.tex  ----------------------------}

\title{On the classifying space of the family of virtually cyclic subgroups}
\author{Wolfgang L\"uck and Michael Weiermann\thanks{\noindent email:
lueck@math.uni-muenster.de,  michi@uni-muenster.de\protect\\
www:~http://www.math.uni-muenster.de/u/lueck/\protect\\
FAX: 49 251 8338370\protect}\\
Fachbereich Mathematik\\ Universit\"at M\"unster\\
Einsteinstr.~62\\ 48149 M\"unster\\Germany}
\maketitle


\typeout{-----------------------  Abstract  ------------------------}

\begin{abstract}
We study the minimal dimension of the classifying space of the family of
virtually cyclic subgroups of a discrete group. We give a complete answer for
instance if the group is virtually poly-$\IZ$, word-hyperbolic or countable
locally virtually cyclic. We give examples of groups for which the difference
of the minimal dimensions of the classifying spaces of virtually cyclic and of
finite subgroups is $-1$, $0$ and $1$, and show in many cases that no other
values can occur.
\\[2mm]
\end{abstract}


 \typeout{--------------------   Section 0: Introduction --------------------------}

\setcounter{section}{-1}
\section{Introduction}

In this paper, we study the classifying spaces $\edub{G} = \EGF{G}{\calvcyc}$
of the family of virtually cyclic subgroups of a group $G$. We are mainly
interested in the minimal dimension $\hdim^G(\edub{G})$ that a $G$-$CW$-model
for $\edub{G}$ can have. The classifying space for proper $G$-actions $\eub{G}
= \EGF{G}{\calfin}$ has already been studied intensively in the literature. The
spaces $\eub{G}$ and $\edub{G}$ appear in the source of the assembly maps whose
bijectivity is predicted by the Baum-Connes conjecture and the Farrell-Jones
conjecture, respectively. Hence the analysis of $\eub{G}$ and $\edub{G}$ is
important for computations of $K$- and $L$-groups of reduced group
$C^*$-algebras and of group rings (see
Section~\ref{sec:Equivariant_homology_of_relative_assembly_maps}). These spaces
can also be viewed as invariants of $G$, and one would like to understand what
properties of $G$ are reflected in the geometry and homotopy theoretic
properties of $\eub{G}$ and $\edub{G}$.

The main results of this manuscript aim at the computation of the homotopy dimension
$\hdim^G(\edub{G})$ (see Definition~\ref{def:homotopy_dimension}), 
which is much more difficult than in the case of
$\hdim^G(\eub{G})$. It has lead to some  (at least for us) surprising
phenomenons, already for some nice groups which can be described explicitly and
have interesting geometry.

After briefly recalling the notion of $\eub{G}$ and $\edub{G}$ in
Section~\ref{sec:Classifying_Spaces_for_Families}, we investigate in
Section~\ref{sec:Passing_to_larger_families} how one can build $\EGF{G}{\calg}$
out of $\EGF{G}{\calf}$ for families $\calf \subseteq \calg$. We construct a
$G$-pushout
\begin{equation*}
\xymatrix{\coprod_{V\in\calm}G\times_{N_GV}\eub{N_GV}
  \ar[d]^{\coprod_{V\in\calm} \id_G\times f_V} \ar[r]^-i & \eub{G} \ar[d] \\
  \coprod_{V\in\calm}G\times_{N_GV}EW_GV \ar[r] & \edub{G}
  }
\end{equation*}
for a certain set $\calm$ of maximal virtually cyclic subgroups of $G$,
provided that $G$ satisfies the condition $(M_{\calfin \subseteq \calvcyc})$
that every infinite virtually cyclic subgroup $H$ is contained in a unique
maximal infinite virtually cyclic subgroup $H_{\max}$.

For a subgroup $H \subseteq G$, we denote by $N_GH:=\{g \in G \mid gHg^{-1} =
H\}$ its \emph{normalizer} and put $W_GH := N_GH/H$. We will say that $G$
satisfies $(N\!M_{\calfin \subseteq \calvcyc})$ if it satisfies $(M_{\calfin
\subseteq \calvcyc})$ and $N_GV = V$, i.e., $W_GV = \{1\}$, holds for all $V
\in \calm$. In
Section~\ref{sec:A_Class_of_Groups_satisfying_(NM_calfin_subseteq_calvcyc)}, we
give a criterion for $G$ to satisfy $(N\!M_{\calfin \subseteq \calvcyc})$.

Section~\ref{sec:Models_for_EGF_for_colimits_of_groups} is devoted to the
construction of models for $\edub{G}$ and $\eub{G}$ from models for
$\edub{G_i}$ and $\eub{G_i}$, respectively, if $G$ is the directed union of the
directed system $\{G_i \mid i \in I\}$ of subgroups.

In Section~\ref{sec:On_the_dimension_of_edub(G)}, we prove various results
concerning $\hdim^G(\edub{G})$ and give some examples. In
Subsection~\ref{subsec:Lower_bounds_for_the_homotopy_dimension_for_calvcyc_by_the_one_for_calfin},
we prove
\begin{equation*}
 \hdim^G(\eub{G}) \le 1 + \hdim^G(\edub{G}).
\end{equation*}
While one knows that $\hdim^{G\times H}(\eub{(G\times H)}) \le
\hdim^{G}(\eub{G}) + \hdim^{H}(\eub{H})$, we prove that
\begin{equation*}
\hdim^{G\times H}(\edub{(G\times H)}) \le \hdim^{G}(\edub{G}) +
\hdim^{H}(\edub{H}) + 3
\end{equation*}
and give examples showing that this inequality cannot be improved in general.

In
Subsection~\ref{subsec:Upper_bounds_for_the_homotopy_dimension_for_calvcyc_by_the_one_for_calfin},
we show
\begin{equation*}
    \hdim^G(\edub{G})
    \begin{cases}
    = \hdim^G(\eub{G}) & \text{if } \hdim^G(\eub{G}) \ge 2; \\
    \le 2 & \text{if } \hdim^G(\eub{G}) \le 1
    \end{cases}
\end{equation*}
provided that $G$ satisfies $(N\!M_{\calfin\subseteq\calvcyc})$, and, more
generally,
\begin{equation*}
  \hdim^G(\edub{G}) \le \hdim^G(\eub{G}) +1
\end{equation*}
provided that $G$ satisfies $(M_{\calfin\subseteq\calvcyc})$. This implies that
the fundamental group $\pi$ of a closed hyperbolic closed manifold $M$
satisfies
\begin{equation*}
    \vcd(\pi) = \dim(N) = \hdim^{\pi}(\eub{\pi}) = \hdim^{\pi}(\edub{\pi}),
\end{equation*}
yielding a counterexample to a conjecture due to Connolly-Fehrmann-Hartglass.

In Subsection~\ref{subsec:Virtually_poly-Z_groups}, a complete computation of
$\hdim^G(\edub{G})$ is presented for virtually poly-$\IZ$ groups. This leads to
some interesting examples in
Subsection~\ref{subsec:Extensions_of_free_abelian_groups}. For instance, we
construct, for $k = -1,0,1$, automorphisms $f_k \colon \Hei \to \Hei$ of the
three-dimensional Heisenberg group $\Hei$ such that
\begin{equation*}
    \hdim^{\Hei\rtimes _{f_k} \IZ}(\edub{(\Hei\rtimes _{f_k} \IZ)}) = 4 + k.
\end{equation*}
Notice that $\hdim^{\Hei\rtimes _{f_k} \IZ}(\eub{(\Hei \rtimes_{f} \IZ)}) =
\cd(\Hei \rtimes_{f} \IZ) = 4$ holds for any automorphism $f \colon \Hei \to
\Hei$.

In Subsection~\ref{subsec:Models_for_EGcalsfg}, we briefly investigate
$\hdim^G(\EGF{G}{\calsfg})$, where $\calsfg$ is the family of subgroups of $G$
which are contained in some finitely generated subgroup. After that, in
Subsection~\ref{subsec:low_dimensions}, we deal with groups for which the
values of $\hdim^G(\eub{G})$ and $\hdim^G(\edub{G})$ are small. We show that a
countable infinite group $G$ is locally finite if and only if $\hdim^G(\eub{G})
= \hdim^G(\edub{G}) = 1$.

Finally, we explain in
Section~\ref{sec:Equivariant_homology_of_relative_assembly_maps} how the models
we construct for $\edub{G}$ can be used to obtain information about the
relative equivariant homology group that appears in the source of the assembly
map in the Farrell-Jones conjecture for algebraic $K$-theory.

All the results in this paper raise the following question:
\begin{quote} \em
For which groups $G$ is it true that
\begin{equation*}
    \hdim^{G}(\eub{G}) - 1 \le \hdim^{G}(\edub{G}) \le \hdim^{G}(\eub{G}) +1 \quad {?}
\end{equation*}
\end{quote}
We have no example of a group for which the above inequality is not true.

In this paper, we are not dealing with the question whether there is a finite
$G$-$CW$-model for $\edub{G}$. Namely, there is no known counterexample to the
conjecture (see~\cite[Conjecture 1]{Juan-Pineda-Leary(2006)}) that a group $G$
possessing a finite $G$-$CW$-model for $\edub{G}$ is virtually cyclic.

The authors wishes to thank the referee for his detailed and valuable report.

The paper was supported by the Sonder\-forschungs\-be\-reich
478 -- Geometrische Strukturen in der Mathematik -- the
Max-Planck-Forschungspreis and the Leibniz-Preis of the first author. Parts of this paper have
already appeared in the Ph.D. thesis of the second author.


\typeout{--------------------   Section 1:   --------------------------}

\section{Classifying Spaces for Families}
\label{sec:Classifying_Spaces_for_Families}

We briefly recall the notions of a family of subgroups and the associated
classifying spaces. For more information, we refer for instance to the survey
article~\cite{Lueck(2005s)}.

A \emph{family $\calf$ of subgroups} of $G$ is a set of  subgroups of $G$ which is
closed under conjugation and taking subgroups. Examples for $\calf$ are
\begin{align*}
    \caltr &= \{\text{trivial subgroup}\}; \\
    \calfin &= \{\text{finite subgroups}\}; \\
    \calvcyc &= \{\text{virtually cyclic subgroups}\}; \\
    \calsfg &= \{\text{subgroups of finitely generated subgroups}\}; \\
    \calall &= \{\text{all subgroups}\}.
\end{align*}

Let $\calf$ be a family of subgroups of $G$. A model for the \emph{classifying
space $\EGF{G}{\calf}$ of the family $\calf$} is a $G$-$CW$-complex $X$ all of
whose isotropy groups belong to $\calf$ and such that for any $G$-$CW$-complex
$Y$ whose isotropy groups belong to $\calf$ there is precisely one $G$-map $Y
\to X$ up to $G$-homotopy. In other words, $X$ is a terminal object in the
$G$-homotopy category of $G$-$CW$-complexes whose isotropy groups belong to
$\calf$. In particular, two models for $\EGF{G}{\calf}$ are $G$-homotopy
equivalent, and for two families $\calf_0 \subseteq \calf_1$ there is precisely
one $G$-map $\EGF{G}{\calf_0} \to \EGF{G}{\calf_1}$ up to $G$-homotopy. There
exists a model for $\EGF{G}{\calf}$ for any group $G$ and any family $\calf$ of
subgroups.

A $G$-$CW$-complex $X$ is a model for $\EGF{G}{\calf}$ if and only if the
$H$-fixed point set $X^H$ is contractible for $H \in \calf$ and is empty for $H
\not\in \calf$.

We abbreviate $\eub{G} := \EGF{G}{\calfin}$ and call it the \emph{universal
$G$-$CW$-complex for proper $G$-actions}. We also abbreviate $\edub{G} :=
\EGF{G}{\calvcyc}$.

A model for $\EGF{G}{\calall}$ is $G/G$. A model for $\EGF{G}{\caltr}$ is the
same as a model for $EG$, which denotes the total space of the universal
$G$-principal bundle $EG \to BG$. Our interest in these notes concerns the
spaces $\edub{G}$.


\typeout{-------------------- Section 2 --------------------------}

\section{Passing to larger families}
\label{sec:Passing_to_larger_families}

In this section, we explain in general how one can construct a model for
$\EGF{G}{\calg}$ from $\EGF{G}{\calf}$ if $\calf$ and $\calg$ are families of
subgroups of the group $G$ with $\calf \subseteq \calg$. Let $\sim$ be an
equivalence relation on $\calg\setminus \calf$ with the properties:
\begin{equation} \label{prop_for_sim}
    \begin{aligned}
    &\bullet \; \text{If } H,K\in\calg\setminus\calf \;\text{with } H\subseteq K,
    \; \text{then } H\sim K; \\
    &\bullet \;\text{If } H,K\in\calg\setminus\calf \;\text{and } g\in G,
    \; \text{then } H\sim K\Leftrightarrow gHg^{-1}\sim gKg^{-1}.
    \end{aligned}
\end{equation}
Let $[\calg \setminus \calf]$ be the set of equivalence classes of $\sim$.
Denote by $[H] \in [\calg \setminus \calf]$ the equivalence class of $H \in \calg
\setminus \calf$, and define the subgroup
\begin{equation*}
    N_G [H] := \{g \in G\mid [g^{-1}Hg] =  [H]\}
\end{equation*}
of $G$. Then $N_G[H]$ is the isotropy group of $[H]$ under the $G$-action on
$[\calg \setminus \calf]$ induced by conjugation. Define a family of subgroups
of $N_G[H]$ by
\begin{equation*}
    \calg [H] := \{K \subseteq N_G[H] \mid K \in \calg\setminus \calf, [K] = [H]\}
    \cup (\calf \cap N_G[H]) ,
\end{equation*}
where $\calf \cap N_G[H] $ consists of all the subgroups of $N_G[H]$ belonging
to $\calf$.

\begin{definition}[Equivalence relation on $\calvcyc\setminus\calfin$]
\label{def:for_sim}
If $\calg = \calvcyc$ and $\calf = \calfin$, we will take for the equivalence
relation $\sim$
\begin{equation*}
    V \sim W \Leftrightarrow |V \cap W| = \infty.
\end{equation*}
\end{definition}

\begin{theorem}[Constructing models from given ones]
\label{the:passing_from_calf_to_calg}
Let $\calf \subseteq \calg$ and $\sim$ be as above such that
properties~\eqref{prop_for_sim} hold. Let $I$ be a complete system of
representatives $[H]$ of the $G$-orbits in $[\calg \setminus \calf]$ under the
$G$-action coming from conjugation. Choose arbitrary $N_G[H]$-$CW$-models for
$\EGF{N_G[H]}{\calf \cap N_G[H]}$ and $\EGF{N_G[H]}{\calg[H]}$, and an
arbitrary $G$-$CW$-model for $\EGF{G}{\calf}$. Define a $G$-$CW$-complex $X$ by
the cellular $G$-pushout
\begin{equation*}
    \xymatrix{\coprod_{[H] \in I} G\times_{N_G[H]} \EGF{N_G[H]}{\calf \cap N_G[H]}
            \ar[d]^{\coprod_{[H] \in I}\id_G \times_{N_G[H]}  f_{[H]}} \ar[r]^-i & \EGF{G}{\calf} \ar[d] \\
                    \coprod_{[H] \in I} G\times_{N_G[H]} \EGF{N_G[H]}{\calg[H]}
      \ar[r] & X
    }
\end{equation*}
such that $f_{[H]}$ is a cellular $N_G[H]$-map for every $[H] \in I$ and $i$ is an
inclusion of $G$-$CW$-complexes, or such that every map $f_{[H]}$ is an inclusion of
$N_G[H]$-$CW$-complexes for every $[H] \in I$ and $i$ is a cellular $G$-map.

Then $X$ is a model for $\EGF{G}{\calg}$.
\end{theorem}
\begin{proof}
We have to show that $X^K$ is contractible if $K$ belongs to $\calg$, and that
it is empty, otherwise. For $[H] \in I$, let $s_{[H]} \colon G/N_G[H] \to G$ be
a set-theoretic section of the projection $G \to G/N_G[H]$.  Let $K \subseteq
G$ be a subgroup. Taking $K$-fixed points in the above $G$-pushout yields, up
to homeomorphism, the following pushout
\begin{equation} \label{eq:pushout_model_2}
\raisebox{2.5em}{\xymatrix{\coprod_{[H] \in I} \coprod_{\alpha}
        \EGF{N_G[H]}{\calf \cap N_G[H]}^{s_{[H]}(\alpha)^{-1}Ks_{[H]}(\alpha)}
        \ar[d]_{\coprod_{[H]\in I}\coprod_{\alpha} f_{[H],\alpha}} \ar[r]^-i
        & \EGF{G}{\calf}^K \ar[d] \\
        \coprod_{[H] \in I} \coprod_{\alpha}
        \EGF{N_G[H]}{\calg[H]}^{s_{[H]}(\alpha)^{-1}Ks_{[H]}(\alpha)} \ar[r] & X^K
    }}
\end{equation}
in which one of the maps starting from the left upper corner is a cofibration
and $\alpha$ runs through $\{\alpha\in G/N_G[H] \mid
s_{[H]}(\alpha)^{-1}Ks_{[H]}(\alpha)\subseteq N_G[H]\}$.

Assume first that $K \notin\calg$. Then the entries in the upper row and the
lower left entry of~\eqref{eq:pushout_model_2} are clearly empty.  Hence, in
this case $X^K$ is empty.

If $K\in\calg\setminus\calf$, the entries in the upper row of
\eqref{eq:pushout_model_2} are again empty, so in order to show that $X^K$ is
contractible, we must show that the lower left entry of
\eqref{eq:pushout_model_2} is contractible. Note that the space
$\EGF{N_G[H]}{\calg[H]}^{s_{[H]}(\alpha)^{-1}Ks_{[H]}(\alpha)}$ is non-empty if
and only if $s_{[H]}(\alpha)^{-1}Ks_{[H]}(\alpha)$ belongs to the family
$\calg[H]$ of subgroups of $N_G[H]$, which is equivalent to the condition that
$[s_{[H]}(\alpha)^{-1}Ks_{[H]}(\alpha)] = [H]$ and
$s_{[H]}(\alpha)^{-1}Ks_{[H]}(\alpha) \subseteq N_G[H]$. Now, there is
precisely one $[H_K] \in I$ for which there exists $g_K \in G$ with the
property that $[g_K^{-1}Kg_K] = [H_K]$. An element $g \in G$ satisfies
$[g^{-1}Kg] = [H_K]$ if and only if $g_K^{-1}g \in N_G[H_K]$. Moreover,
$L\subseteq N_G[H]$ always holds if $L\in\calg\setminus\calf$ is such that
$[L]=[H]$. We conclude from these observations that
$\EGF{N_G[H]}{\calg[H]}^{s_{[H]}(\alpha)^{-1}Ks_{[H]}(\alpha)}$ is non-empty if
and only if $[H] = [H_K]$ and $\alpha = g_KN_G[H_K]$. Furthermore, it is
contractible in this case.  This implies that the lower left entry of
\eqref{eq:pushout_model_2} is contractible.

Finally, if $K\in\calf$, then the upper right entry of
\eqref{eq:pushout_model_2} is contractible. This implies that the same will
hold for $X^K$ if we can show that all the maps $f_{[H],\alpha}$ are homotopy
equivalences. But this is clear since the source and target spaces of the
$f_{[H],\alpha}$ are contractible.
\end{proof}

\begin{remark}[Dimension of the model constructed for $\EGF{G}{\calg}$]
\label{rem:how_to_apply_theorem_the:passing_from_calf_to_calg}
A $G$-pushout as described in Theorem~\ref{the:passing_from_calf_to_calg}
always exists. To be more precise, the maps $i$ and $f_{[H]}$ exist and are
unique up to equivariant homotopy equivalence due to the universal property of
the classifying spaces $\EGF{G}{\calf}$ and $\EGF{N_G[H]}{\calg[H]}$. Moreover,
because of the equivariant cellular approximation theorem (see for
instance~\cite[Theorem~2.1 on page~32]{Lueck(1989)}), these maps can be assumed
to be cellular. Finally, one can replace a cellular $G$-map $f \colon X \to Y$
by the canonical inclusion $j \colon X \to \cyl(f)$ into its mapping cylinder.
Then the canonical projection $p \colon \cyl(f) \to Y$ is a $G$-homotopy
equivalence such that $\pr \circ i = f$.

Note that the dimension of $\cyl(f)$ is $\max\{1 + \dim(X),\dim(Y)\}$. This
means that we we get the following conclusion from
Theorem~\ref{the:passing_from_calf_to_calg} which we will often use: There
exists an $n$-dimensional $G$-$CW$-model for $\EGF{G}{\calg}$ if there exists
an $n$-dimensional $G$-$CW$-model for $\EGF{G}{\calf}$ and, for every $H\in I$,
an $(n-1)$-dimensional $N_G[H]$-$CW$-model for $\EGF{N_G[H]}{\calf \cap
N_G[H]}$ and an $n$-dimensional $N_G[H]$-$CW$-model for
$\EGF{N_G[H]}{\calg[H]}$.
\end{remark}

\begin{example} \label{exa:N_G[V]_is_N_GC}
In general $N_G[V]$, cannot be written as a normalizer $N_GW$ for some $W \in
[V]$ in the case $\calg = \calvcyc$, $\calf = \calfin$ and $\sim$ as defined in
Definition~\ref{def:for_sim}.

For instance, let $p$ be a prime number. Let $\IZ[1/p]$ be the subgroup of
$\IQ$ consisting of rational numbers $x \in \IQ$ for which $p^n \cdot x \in
\IZ$ for some positive integer $n$. This is the directed union $\bigcup_{n \in
\IN} p^{-n} \cdot \IZ$.  Let $p\cdot \id \colon \IZ[1/p] \to \IZ[1/p]$ be the
automorphism given by multiplication with $p$, and define $G$ to be the
semi-direct product $\IZ[1/p] \rtimes_{p\cdot \id} \IZ$. Let $C$ be the cyclic
subgroup of $\IZ[1/p]$ generated by $1$. Then $N_G[C] = G$ but $N_GV =
\IZ[1/p]$ for every $V \in [C]$.
\end{example}

\begin{notation} \label{not:(Mcalg,calf)}
Let $\calf\subseteq\calg$ be families of subgroups of a group $G$.

We shall say that \emph{$G$ satisfies $(M_{\calf \subseteq \calg})$} if every
subgroup $H\in\calg\setminus\calf$ is contained in a unique
$H_{\max}\in\calg\setminus\calf$ which is maximal in $\calg\setminus\calf$,
i.e., $H_{\max} \subseteq K$ for $K\in\calg\setminus\calf$ implies
$H_{\max}=K$.

We shall say that \emph{$G$ satisfies $(N\!M_{\calf \subseteq \calg})$} if $M$
satisfies $(M_{\calf \subseteq \calg})$ and every maximal element $M \in
\calg\setminus \calf$ equals its normalizer, i.e., $N_GM = M$ or, equivalently,
$W_GM=\{1\}$.
\end{notation}

Examples of groups satisfying $(N\!M_{\calfin \subseteq \calvcyc})$ will be
given in
Section~\ref{sec:A_Class_of_Groups_satisfying_(NM_calfin_subseteq_calvcyc)}. In
\cite[page 101]{Davis-Lueck(2003)}, Davis-L\"uck discuss examples of groups
satisfying $(M_{\caltr \subseteq \calfin})$ are discussed. The following will
be a consequence of Theorem~\ref{the:passing_from_calf_to_calg}:

\begin{corollary} \label{cor:passing_from_calf_to_calg_assuming_(McalG,calF)}
Let $\calf\subseteq\calg$ be families of subgroups of a group $G$ which
satisfies $(M_{\calf \subseteq \calg})$. Let $\calm$ be a complete system of
representatives of the conjugacy classes of subgroups in $\calg \setminus
\calf$ which are maximal in $\calg \setminus \calf$.  Let $\calsub(M)$ be the
family of subgroups of $M$. Consider a cellular $G$-pushout
\begin{equation*}
    \xymatrix{\coprod_{M \in \calm} G\times_{N_GM} \EGF{N_GM}{\calf \cap N_GM}
        \ar[d]^{\coprod_{M \in \calm} \id_G \times_{N_GM}  f_{[H]}} \ar[r]^-i
        & \EGF{G}{\calf} \ar[d] \\
        \coprod_{M \in \calm} G\times_{N_GM} \EGF{N_GM}{\calsub(M) \cup
        (\calf \cap N_GM)} \ar[r] & X}
\end{equation*}
such that $f_{[H]}$ is a cellular $N_G[H]$-map for every $[H] \in I$ and $i$ is
an inclusion of $G$-$CW$-complexes, or such that $f_{[H]}$ is an inclusion of
$N_G[H]$-$CW$-complexes for every $[H] \in I$ and $i$ is a cellular $G$-map.

Then $X$ is a model for $\EGF{G}{\calg}$.
\end{corollary}

\begin{proof}
Let $H_{\max}$ denote the unique maximal element in $\calg\setminus\calf$ which
contains $H\in\calg\setminus\calf$. We use the equivalence relation on
$\calg\setminus \calf$ given by $H \sim K \Leftrightarrow H_{\max} = K_{\max}$.
Then, for $M \in \calm$, the family $\calg[M]$ is just $\{K \in \calg \mid K
\subseteq M \text{ or } K \in \calf\}$. Now apply
Theorem~\ref{the:passing_from_calf_to_calg}.
\end{proof}

\begin{remark}
If $G$ satisfies $(M_{\calfin\subseteq \calvcyc})$, then the equivalence
relation $\sim$ defined in Definition~\ref{def:for_sim} agrees with the
equivalence relation $H \sim K \Leftrightarrow H_{\max} = K_{\max}$ appearing
in the proof of
Corollary~\ref{cor:passing_from_calf_to_calg_assuming_(McalG,calF)} above.
\end{remark}

\begin{corollary} \label{cor:Efin_and_Evcyc_model_assuming_(McalG,calF)}
Let $G$ be a group satisfying $(M_{\caltr\subseteq \calfin})$ or
$(M_{\calfin\subseteq \calvcyc})$ respectively. We denote by $\calm$ a complete system of
representatives of the conjugacy classes of maximal finite subgroups
$F\subseteq G$ or of maximal infinite virtually cyclic subgroups $V\subseteq
G$ respectively. Consider the cellular $G$-pushouts
\begin{equation*}
  \xymatrix{\coprod_{F\in\calm} G\times_{N_GF}EN_GF
    \ar[d]^{\coprod_{F\in\calm}\id_G \times f_F} \ar[r]^-i & EG \ar[d] \\
    \coprod_{F\in\calm} G\times_{N_GF}EW_GF \ar[r] & X }
  \raisebox{-2.0em}{\quad\;\text{or}\quad}
  \xymatrix{\coprod_{V\in\calm}G\times_{N_GV}\eub{N_GV}
    \ar[d]^{\coprod_{V\in\calm} \id_G\times f_V} \ar[r]^-i & \eub{G} \ar[d] \\
    \coprod_{V\in\calm}G\times_{N_GV}EW_GV \ar[r] & Y }
\end{equation*}
where $EW_GH$ is viewed as an $N_GH$-$CW$-complex by restricting with the
projection $N_GH \to W_GH$ for $H \subseteq G$, the maps starting from the left
upper corner are cellular and one of them is an inclusion of
$G$-$CW$-complexes.

Then $X$ is a model for $\eub{G}$ or $Y$ is a model for $\edub{G}$ respectively.
\end{corollary}
\begin{proof}
This follows from
Corollary~\ref{cor:passing_from_calf_to_calg_assuming_(McalG,calF)} and the
following facts. A model for $\EGF{N_GH}{\calsub(H)}$ is given by $EW_GH$
considered as an $N_GH$-$CW$-complex by restricting with the projection $N_GH
\to W_GH$ for every $H \subseteq G$. Furthermore, we obviously have $\caltr
\cap N_GF \subseteq \calsub(F)$ for every maximal finite subgroup $F$, whereas
$\calfin \cap N_GV \subseteq \calsub(V)$ for every maximal infinite virtually
cyclic subgroup $V\subseteq G$ because $W_GV$ is torsionfree (if it were not,
then the preimage of a non-trivial finite subgroup of $W_GV$ under the
projection $N_GV\to W_GV$ would be virtually cyclic and strictly containing
$V$, contradicting its maximality).
\end{proof}

As a special case of
Corollary~\ref{cor:Efin_and_Evcyc_model_assuming_(McalG,calF)}, we get

\begin{corollary} \label{cor:Efin_and_Evcyc_model_with_N_GM_is_M}
Let $G$ be a group satisfying $(N\!M_{\caltr\subseteq \calfin})$ or
$(N\!M_{\calfin\subseteq \calvcyc})$ respectively. Let $\calm$ be a complete system of
representatives of the conjugacy classes of maximal finite subgroups
$F\subseteq G$ or of maximal infinite virtually cyclic subgroups $V\subseteq
G$ respectively. Consider the cellular $G$-pushouts
\begin{equation*}
    \xymatrix{\coprod_{F\in\calm} G\times_{F}EF
        \ar[d]^{\coprod_{F\in\calm} p} \ar[r]^-i & EG \ar[d] \\
        \coprod_{F\in\calm} G/F \ar[r] & X
    }
    \raisebox{-2.0em}{\quad\;\text{or}\quad}
    \xymatrix{\coprod_{V\in\calm}G\times_{V}\eub{V}
    \ar[d]^{\coprod_{V\in\calm} p} \ar[r]^-i & \eub{G} \ar[d] \\
    \coprod_{V\in\calm}G/V \ar[r] & Y
    }
\end{equation*}
where $i$ is an inclusion of $G$-$CW$-complexes and $p$ is the obvious
projection.

Then $X$ is a model for $\eub{G}$ or  $Y$ is a model for $\edub{G}$ respectively.
\end{corollary}


\typeout{--------------------   Section 3   --------------------------}

\section{A Class of Groups satisfying $(N\!M_{\calfin \subseteq \calvcyc})$}
\label{sec:A_Class_of_Groups_satisfying_(NM_calfin_subseteq_calvcyc)}

The following provides a criterion for a group to satisfy
$(N\!M_{\calf\subseteq \calvcyc})$:
\begin{theorem} \label{the:max_vc_sg}
Suppose that the group $G$ satisfies the following two conditions:
\begin{itemize}
\item Every infinite cyclic subgroup $C\subseteq G$ has finite index
  $[C_GC:C]$ in its centralizer;
\item Every ascending chain $H_1\subseteq H_2\subseteq\ldots$ of finite
  subgroups of $G$ becomes stationary, i.e., there is an $n_0\in \IN$
  such that $H_n=H_{n_0}$ for all $n\geq n_0$.
\end{itemize}
Then every infinite virtually cyclic subgroup $V\subseteq G$ is contained in a
unique maximal virtually cyclic subgroup $V_{\max}\subseteq G$. Moreover,
$V_{\max}$ is equal to its normalizer $N_G(V_{\max})$, and
\begin{equation*}
    V_{\max} = \bigcup_{C\subseteq V}N_GC,
\end{equation*}
where the union is over all infinite cyclic subgroups $C\subseteq V$ that are
normal in $V$.

In particular, $G$ satisfies $(N\!M_{\calfin\subseteq \calvcyc})$.
\end{theorem}
\begin{proof}
We fix an infinite virtually cyclic subgroup $V\subseteq G$ and define
$V_{\max} := \bigcup_{C\subseteq V}N_GC$, where the union is over all infinite
cyclic subgroups $C\subseteq V$ that are normal in $V$. The collection of all
such subgroups of $V$ is countable, and we denote it by $\{C_n\}_{n\in\IN}$.
Since every index $[V:C_n]$ is finite, $[V:C_0\cap\ldots\cap C_n]$ must also be
finite for $n\in\IN$. Thus, if we set $Z_n:= C_0\cap\ldots\cap C_n$, then
$Z_n\subseteq C_n$, and $Z_0\supseteq Z_1\supseteq\ldots$ is a descending chain
of infinite cyclic subgroups of $V$ that are normal in $V$. If $C'\subseteq C$
are two infinite cyclic subgroups of $G$, then $N_GC\subseteq N_GC'$ because
$C'$ is a characteristic subgroup of $C$. It follows in our situation that
$N_GC_n\subseteq N_GZ_n$, which implies
\begin{equation} \label{eq:max_vc_sg}
    \bigcup_{n=0}^\infty N_GC_n = \bigcup_{n=0}^\infty N_GZ_n.
\end{equation}
Furthermore, $N_GZ_0\subseteq N_GZ_1\subseteq\ldots$ is an ascending chain,
which becomes stationary by the following argument. Namely, we can estimate
\begin{equation*}
  [N_GZ_n:N_GZ_0] \le [N_GZ_n:C_GZ_0] = [N_GZ_n:C_GZ_n]\cdot[C_GZ_n:C_GZ_0],
\end{equation*}
and the first factor on the right is not greater than $2$ since there is an
injection $N_GZ_n/C_GZ_n \to \aut(Z_n)$, whereas the second does not exceed an
appropriate constant as we will show in Lemma~\ref{lem:upper_bound} 
below. Hence, it follows from
\eqref{eq:max_vc_sg} that $V_{\max}=N_GZ_n$ for all $n\geq n_0$ if $n_0\in\IN$
is sufficiently large. In other words, we can record that this construction
yields, for any infinite cyclic subgroup $C\subseteq V$ that is normal in $V$,
an infinite cyclic subgroup $Z\subseteq C$ that is normal in $V$ such that for
all infinite cyclic subgroups $Z'\subseteq Z$ we have $V_{\max}=N_GZ'$.

Using the equality $V_{\max}=N_GZ_{n_0}$, it is now obvious that $V_{\max}$ is
virtually cyclic since the finite index subgroup $C_GZ_{n_0}$ is so due to the
assumption imposed on $G$. In addition, $V\subseteq V_{\max}$ since $Z_{n_0}$
is normal in $V$. Now suppose $W\subseteq G$ is an infinite virtually cyclic
subgroup such that $V\subseteq W$. We claim that $V_{\max}=W_{\max}$. In order
to prove this, let $Z_W\subseteq W$ be an infinite cyclic subgroup that is
normal in $W$ such that $W_{\max}=N_GZ_W'$ holds for all infinite cyclic
subgroups $Z_W'\subseteq Z_W$.  Then $Z_W\cap V$ is a finite index subgroup of
$W$ and so is an infinite cyclic subgroup of $V$ that is normal in $V$. As we
have already seen, there is an infinite cyclic subgroup $Z_V\subseteq Z_W\cap
V$ that is normal in $V$ such that $V_{\max}=N_GZ_V'$ holds for all infinite
cyclic subgroups $Z_V'\subseteq Z_V$.  Since $V$ has finite index in $W$ and
$Z_V$ is normal in $V$, the intersection $Z_V''$ of all the conjugates of $Z_V$
by elements in $W$ has finite index in $Z_V$ and therefore is an infinite
cyclic subgroup of $Z_V$.  This implies $V_{\max} = N_GZ_V'' = W_{\max}$. From
this we can deduce immediately that $V_{\max}$ is indeed maximal among
virtually cyclic subgroups of $G$ containing $V$ and that it is uniquely
determined by this property.

Finally, we will show that $N_GV_{\max}$ is virtually cyclic, so that it is
equal to $V_{\max}$. Let $C\subseteq V_{\max}$ be infinite cyclic. Since $C$
has finite index in $V_{\max}$ and $V_{\max}$ contains only finitely many
subgroups of index $[V_{\max}:C]$, the group $D$ which we define as the
intersection of all conjugates of $C$ in $N_GV_{\max}$ has finite index in
$V_{\max}$ and is therefore infinite cyclic as well. Obviously, $D$ is normal
in $N_GV_{\max}$, so $N_GV_{\max}\subseteq N_GD$ holds, the latter being
virtually cyclic because the finite index subgroup $C_GD$ is so by assumption.
Hence Theorem~\ref{the:max_vc_sg} is proven as soon as we have finished the proof
of the next lemma.
\end{proof}

\begin{lemma} \label{lem:upper_bound}
There exists a natural number $N$ such that for all $n \in \IN$ we have
$$[C_GZ_n:C_GZ_0] \le N.$$
\end{lemma}
\begin{proof}
Since the center of the virtually cyclic
group $C_GZ_n$ contains $Z_n$ and hence is infinite, $C_GZ_n$ is of type I,
i.e., possesses an epimorphism onto $\IZ$. Denoting by $T_n$ the torsion
subgroup of $H_1(C_GZ_n)$, it then follows that $H_1(C_GZ_n)/T_n$ is infinite
cyclic. Let $p_n\colon C_GZ_n\to H_1(C_GZ_n)\to H_1(C_GZ_n)/T_n$ be the
canonical projection, and let $c_n\in C_GZ_n$ be such that $p_n(c_n)$ is a
generator.

Consider the commutative diagram
\begin{equation*}
    \xymatrix{1 \ar[r] & \ker(p_0) \ar[d] \ar[r] & C_GZ_0 \ar[d] \ar[r]^-{p_0} &
            H_1\bigl(C_GZ_0\bigr)/T_0 \ar[d] \ar[r] & 1 \\
        1 \ar[r] & \ker(p_1) \ar[d] \ar[r] & C_GZ_1 \ar[d] \ar[r]^-{p_1} &
            H_1\bigl(C_GZ_1\bigr)/T_1 \ar[d] \ar[r] & 1 \\
        & \vdots & \vdots & \vdots
    }
\end{equation*}
which has exact rows and in which all the vertical arrows in the first and
second column are inclusions and all vertical arrows in the third column are
induced by the obvious inclusions.  They, too, are injective because
$\ker(p_n)$ is finite and $H_1(C_GZ_n)/T_n$ is infinite cyclic for $n\in\IN$.
Since all the $\ker(p_n)$ are finite, $\bigcup_{n=0}^\infty\ker(p_n)$ is again
finite, say of order $a\in\IN$, by the assumption on $G$. Then in particular
$[\ker(p_n):\ker(p_0)]\le a$ for $n\in\IN$, so
$\bigl\{[C_GZ_n:C_GZ_0]\,\big|\,n\in\IN\bigl\}$ will be bounded if we can show
that the index of $H_1(C_GZ_0)/T_0$ in $H_1(C_GZ_n)/T_n$ is less than a
constant which does not depend on $n\in\IN$.

In order to construct such a constant, let $r_n\in\IZ$ be such that $p_0(c_0)$
is mapped to $p_n(c_n)^{r_n}$ under the inclusions in the above diagram. Then,
by exactness, there is a $k_n\in\ker(p_n)$ such that $c_0=k_nc_n^{r_n}$. Since
the order of $\ker(p_n)$ divides $a$, the group $\aut\bigl(\ker(p_n)\bigr)$
contains at most ${a!}$  elements, so that any
$\phi_n\in\aut\bigl(\ker(p_n)\bigr)$ satisfies $\phi_n^{a!}=\id$. Put
$d:=a\cdot {a!}$, then we get for every $\phi_n\in\aut\bigl(\ker(p_n)\bigr)$
and $k\in\ker(p_n)$ the equality
\begin{equation*}
    \prod_{i=0}^{d-1}\phi_n^i(k)
    = \prod_{j=1}^a\,\prod_{i=(j-1)\cdot a!}^{j\cdot a!-1}\phi_n^i(k)
    = \left(\prod_{i=0}^{a!-1}\phi_n^i(k)\right)^a = 1.
\end{equation*}
Specializing to $\phi_n(k):= c_n^{r_n}kc_n^{-r_n}$ yields
\begin{equation*}
    c_0^d = \left(k_n c_n^{r_n}\right)^d
    = \left(\prod_{i=0}^{d-1}\phi_n^i(k_n)\right)\cdot c_n^{r_nd}=c_n^{r_nd}.
\end{equation*}
If $z_0$ is a generator of $Z_0$, then there exists an $s\in\IZ$ such that
$p_0(z_0)=p_0(c_0)^s$. This means $z_0^{-1}c_0^s\in\ker(p_0)$, hence
$z_0^d=c_0^{sd}$. Altogether this implies that if $dZ_0$ denotes the cyclic
group generated by $z_0^d$, we have $dZ_0=\langle c_n^{r_nds}\rangle$ and thus
$c_n\in C_GdZ_0$.

We can finally define $b:=[C_GdZ_0:dZ_0]$, which is finite by assumption and
constitutes the required constant. This is due to the fact that $c_n^b\in
dZ_0$, so there is a $t_n\in\IZ$ such that $c_n^b=c_n^{r_n dst_n}$. Hence $r_n$
divides $b$, and the index of $H_1(C_GZ_0)/T_0$ in $H_1(C_GZ_n)/T_n$ equals
$|r_n|$ by construction.
\end{proof}

\begin{lemma}\label{lem:pin_leary}
Let $G$ be a group with the property that every non-virtually cyclic subgroup
of $G$ contains a copy of $\IZ\ast\IZ$. Then $G$ satisfies the conditions
appearing in Theorem~\ref{the:max_vc_sg}.
\end{lemma}
\begin{proof}
Let $C\subseteq G$ be infinite cyclic. By assumption, its centralizer $C_GC$ is
either virtually cyclic or contains $\IZ\ast\IZ$. To prove the first condition
appearing in Theorem~\ref{the:max_vc_sg}, we have to show that $C_GC$ does not
contain $\IZ\ast\IZ$. Suppose it does.  Then $\IZ\ast\IZ\cap C=\{1\}$ as
$\IZ\ast\IZ$ has a trivial center. Hence one of the generators of $\IZ\ast\IZ$
together with a generator of $C$ generate a copy of $\IZ\oplus\IZ$ inside $G$,
which contradicts the assumption imposed on $G$.

It is obvious that any ascending chain $H_1\subseteq H_2\subseteq\ldots$ of
finite subgroups of $G$ must become stationary since, otherwise, $\bigcup_nH_n$
would be an infinite torsion subgroup of $G$, contradicting the assumptions on
$G$. This proves the second condition appearing in Theorem~\ref{the:max_vc_sg}.
\end{proof}

\begin{remark}
Lemma~\ref{lem:pin_leary} implies that the main result of Juan-Pineda and Leary
in \cite{Juan-Pineda-Leary(2006)} is covered by Theorem~\ref{the:max_vc_sg}.
Obviously this is also true for \cite[Theorem~8.11]{Lueck(2005s)} which
motivated Theorem~\ref{the:max_vc_sg}.

It follows from the Kurosh subgroup theorem (see for
example~\cite[Theorem~I.5.14 on page~56]{Serre(1980)}) that the class of groups
satisfying the conditions of Theorem~\ref{the:max_vc_sg} is closed under
arbitrary free products, whereas this is not the case for the class of groups
considered in \cite[Theorem~8.11]{Lueck(2005s)}.
\end{remark}

\begin{example}[Word-hyperbolic groups] \label{exa:word-hyp}
Each word-hyperbolic group $G$ satisfies the two conditions appearing in
Theorem~\ref{the:max_vc_sg}. A proof for this can be found in \cite[Theorem~3.2
in~III.$\Gamma$.3 on page~459 and Corollary~3.10 in III.$\Gamma$.3 on
page~462]{Bridson-Haefliger(1999)}.
\end{example}

\begin{example}[A group not satisfying $(M_{\calfin\subseteq \calvcyc})$]
 \label{exa:groups_not_satisfying_M(calfin_subseteq_calvcyc)}
Let $G$ be the semi-direct product $\IZ \rtimes \IZ$ with respect to the group
homomorphism $\IZ \to \aut(\IZ)$ sending the generator to $-\id \colon \IZ \to
\IZ$. A presentation of $G$ is given by $\langle s,t \mid sts^{-1}=
t^{-1}\rangle$. Consider the infinite cyclic subgroups $\langle s \rangle$ and
$\langle st \rangle$. They are maximal among infinite virtually cyclic
subgroups. These two groups are obviously different. Their intersection is the
proper subgroup $\langle s^2 \rangle$. Hence $G$ does \emph{not} satisfy
$(M_{\calfin\subseteq \calvcyc})$.
\end{example}


\typeout{-------------------- Section 4 --------------------------}

\section{Models for $\EGF{G}{\calf}$ for colimits of groups}
\label{sec:Models_for_EGF_for_colimits_of_groups}

In this section, we want to construct models for $\EGF{G}{\calf}$ if $G$ is a
colimit of a directed system of groups $\{G_i \mid i \in I\}$.

Given a directed set $I = (I,\le)$, we consider $I$ as a category. Namely, the
set of objects is the set $I$ itself and $\mor_I(i,j)$ consists of precisely
one element if $i \le j$ and is empty otherwise. A directed system of groups is
a covariant functor from the category $I$ into the category of groups, i.e., it
consists of a collection $\{G_i \mid i \in I\}$ of groups together with group
homomorphisms $\phi_{i,j}\colon G_i\to G_j$ for $i\le j$. Its colimit $G:=
\colim_{i \in I} G_i$ is a group $G$ together with group homomorphisms $\psi_i
\colon G_i \to G$ for every $i \in I$ which satisfy $\psi_j = \psi_i \circ
\phi_{i,j}$ for $i,j \in I$ and $i \le j$ such that the following universal
property holds: Given a group $K$ together with group homomorphisms $\mu_i
\colon G_i \to K$ satisfying $\mu_j = \mu_i \circ \phi_{i,j}$ for $i,j \in I$
and $i \le j$, there is precisely one group homomorphism $\nu \colon G \to K$
satisfying $\nu \circ \psi_i = \mu_i$ for every $i \in I$. If $\{G_i \mid i \in
I\}$ is a directed system of subgroups of $G$ directed by inclusion, then $G =
\colim_{i \in I} G_i$ if and only if $G = \bigcup_{i \in I} G_i$.

Given a homomorphism $\psi \colon G' \to G$ of groups and a family of subgroups
$\calf$ of $G$, define $\psi^*\calf$ to be the family of subgroups of $G'$
given by $\{H \subseteq GÄ' \mid \psi(H) \in \calf\}$. If $\psi$ is an
inclusion of groups, then $\psi^*\calf$ agrees with the family $\calf \cap G'$
consisting of all the subgroups of $G'$ which belong to $\calf$.

\begin{definition}[Homotopy dimension] \label{def:homotopy_dimension}
Given a $G$-space $X$, the \emph{homotopy dimension} $\hdim^G(X) \in \{0,1,
\ldots \} \amalg \{\infty\}$ of $X$ is defined to be the infimum over the
dimensions of all $G$-$CW$-complexes $Y$ which are $G$-homotopy equivalent to
$X$.
\end{definition}

Obviously, $\hdim^G(X)$ depends only on the $G$-homotopy type of $X$.  In
particular, $\hdim^G(\EGF{G}{\calf})$ is well-defined for any group $G$ and any
family $\calf$. For example, $\hdim^G(\EGF{G}{\calf}) = 0$ is equivalent to the
condition that $\calf = \calall$.

Let $\calc$ be a small category. A contravariant $\calc$-$CW$-complex is a
contravariant functor from $\calc$ to the category of $CW$-complexes which is
built of cells of the shape $\mor_\calc(?,c)\times D^n$. For a rigorous
definition of this object, we refer for instance
to~\cite[Definition~3.2]{Davis-Lueck(1998)} (where it is called a contravariant
free $\calc$-$CW$-complex). A model for the \emph{contravariant classifying
space $E\calc$ of the category $\calc$} is, by definition, a contravariant
$\calc$-$CW$-complex whose evaluation at any object in $\calc$ is contractible.
Let $\hdim(\calc)$ be the infimum of the dimensions of all models for $E\calc$.
The elementary proof of the next result is left to the reader
(see~\cite[Example~3.9]{Davis-Lueck(1998)})).

\begin{lemma} \label{lem:hdim(I)}
  Let $I$ be a directed set. Then $I$ contains a maximal element if and only if $\hdim(I)
  = 0$. If $I$ contains a countable cofinal subset $I'$, then $\hdim(I) \le 1$. (Cofinal
  means that for every $i \in I$ there exists $i' \in I$ with $i \le i'$.)
\end{lemma}

Now we can formulate the main result of this section, whose proof is carried
out in the remainder of this section.

\begin{theorem} \label{the:model_for_colimits}
  Let $\{G_i \mid i \in I\}$ be a directed system of groups (with not necessarily
  injective structure maps $\phi_{i,j} \colon G_i \to G_j$ for $i \le j$).  Let $G =
  \colim_{i \in I} G_i$ be the colimit of this directed system.  Let $\calf$ be a family
  of subgroups of $G$ such that every $H \in \calf$ is contained in $\im(\psi_i)$ for some
  $i \in I$, where $\psi_i \colon G_i \to G$ is the structure map of the colimit for $i
  \in I$.  Then
  $$\hdim^G(\EGF{G}{\calf}) \le \hdim(I) +
  \sup\big\{\hdim^{G_i}(\EGF{G_i}{\psi_i^*\calf})\big| i \in I\big\}.$$
\end{theorem}

We need some preparations for the proof of
Theorem~\ref{the:model_for_colimits}. In the following, we will use the
notation and results of~\cite[Section~1 and~3]{Davis-Lueck(1998)}, except that
we replace the notation $\otimes_{\calc}$
in~\cite[Definition~1.4]{Davis-Lueck(1998)} by $\times_{\calc}$. We also
mention that, for a covariant functor $X \colon \calc \to \SPACES$, the space
$E\calc \times_{\calc} X$ is a model for the homotopy colimit of $X$
(see~\cite[page~225]{Davis-Lueck(1998)}).

\begin{lemma} \label{lem:Z_times_calc_X}
Let $\calc$ be a small category. Let $Z$ be a contravariant
$\calc$-$CW$-complex and let $X \colon \calc \to G\text{-}\SPACES$ be a
covariant functor from $\calc$ to the category of $G$-spaces.

Then $Z \times_{\calc} X$ is a $G$-space and we have:
\begin{enumerate}

\item \label{lem:Z_times_calc_X:G-CW-Structure}
If $X(c)$ is a $G$-$CW$-complex for every $c \in \calc$, then $Z \times_{\calc}
X$ has the homotopy type of a $G$-$CW$-complex. We have
\begin{equation*}
    \hdim^G(Z\times_{\calc}X)\le\dim(Z)+\sup\big\{\hdim^{G}\bigl(X(c)\bigr)
    \,\big|\, c\in\ob(\calc)\big\};
\end{equation*}

\item \label{lem:Z_times_calc_X:fixed_point_sets}
For every $H \subseteq G$ there is a natural homeomorphism of spaces
\begin{equation*}
    Z \times_{\calc} X^H \xrightarrow{\cong} \left( Z \times_{\calc} X\right)^H.
\end{equation*}

\end{enumerate}
\end{lemma}
\begin{proof}\ref{lem:Z_times_calc_X:G-CW-Structure}
Let $Z_n$ be the $n$-skeleton of $Z_n$. We construct, by induction over $n =
-1,0,1,\ldots$, a sequence of inclusions of $G$-$CW$-complexes $Y_{-1}
\subseteq Y_0 \subseteq Y_1 \subseteq \ldots$ together with $G$-homotopy
equivalences $f_n \colon Z_n \times_{\calc} X \to Y_n$ such that $f_n|_{Z_{n-1}
\times_{\calc} X} = f_{n-1}$ and
\begin{equation*}
    \dim(Y_n) \le n+\sup\big\{\hdim^{G}\bigl(X(c)\bigr) \,\big|\, c\in\ob(\calc)\big\}
\end{equation*}
holds for $n = 0,1,2, \ldots{}$ Then the claim is obviously true if $Z$ is
finite-di\-men\-sional. For infinite-dimensional $Z$, a standard colimit
argument implies that the system $(f_n)_{n \in \IN}$ of $G$-homotopy
equivalences yields a $G$-homotopy equivalence $f \colon Z \times_{\calc} X \to
Y:= \bigcup_{n \in \IN} Y_n$ since every map $Z_{n-1} \times_{\calc} X \to Z_n
\times_{\calc} X$ is a $G$-cofibration and $Z \times_{\calc} X$ is $\bigcup_{n
\in \IN} Z_n \times_{\calc} X_n$ with respect to the weak topology. As $Y$ is a
$G$-$CW$-complex, the claim follows in this case as well.

The induction begins by taking $Y_{-1} = \emptyset$. The induction step from
$n-1$ to $n \ge 0$ is done as follows.  We can write $Z_n$ as a pushout of
contravariant $\calc$-spaces
\begin{equation*}
    \xymatrix{\coprod_{j \in J}\mor_{\calc}(?,c_j)\times S^{n-1}
         \ar[d]_{\coprod_{j\in J}\id_{\mor_{\calc}(?,c_j)}\times i} \ar[r]
         & Z_{n-1} \ar[d] \\
        \coprod_{j \in J}\mor_{\calc}(?,c_j)\times D^n \ar[r] & Z_n
    }
\end{equation*}
where $i \colon S^{n-1} \to D^n$ is the inclusion.  Applying $- \times_{\calc}
X$ yields a $G$-pushout
\begin{equation*}
    \xymatrix{\coprod_{j \in J}X(c_j)\times S^{n-1}
         \ar[d]_{\coprod_{j\in J}\id_{X(c_j)}\times i} \ar[r]^-u
         & Z_{n-1}\times_\calc X \ar[d] \\
        \coprod_{j \in J}X(c_j)\times D^n \ar[r] & Z_n\times_\calc X
    }
\end{equation*}
Obviously, the left vertical arrow is a $G$-cofibration and hence the same is
true for the right vertical arrow. Choose, for every $c \in \calc$, a
$G$-$CW$-complex $X(c)'$ of dimension $\hdim^G(X(c))$ together with a
$G$-homotopy equivalence $u_c \colon X(c)' \to X(c)$.  By the equivariant
cellular approximation theorem (see for instance~\cite[Theorem~2.1 on
  page~32]{Lueck(1989)}), the composition
\begin{equation*}
  \coprod_{j \in J} X(c_j)' \times S^{n-1} \xrightarrow{\coprod_{j \in J} u_{c_j} \times
    \id_{S^{n-1}}} \coprod_{j \in J} X(c_j) \times S^{n-1} \xrightarrow{u} Z_{n-1}
  \times_{\calc} X \xrightarrow{f_{n-1}} Y_{n-1}
\end{equation*}
is $G$-homotopic to some cellular $G$-map $v \colon \coprod_{j \in J} X(c_j)'
\times S^{n-1} \to Y_{n-1}$. Define $Y_n$ by the $G$-pushout
\begin{equation*}
    \xymatrix{\coprod_{j \in J}X(c_j)'\times S^{n-1}
         \ar[d]_{\coprod_{j\in J}\id_{X(c_j)'}\times i} \ar[r]^-v & Y_{n-1} \ar[d] \\
        \coprod_{j \in J}X(c_j)'\times D^n \ar[r] & Y_n
    }
\end{equation*}
This is a $G$-$CW$-complex whose dimension satisfies
\begin{equation*}
    \dim(Y_n) \le \max \big\{n+\sup\big\{\dim\bigl(X(c_j)'\bigr) \,\big|\,
    j\in J\big\}, \dim(Y_{n-1})\big\},
\end{equation*}
and the induction hypothesis applied to $Y_{n-1}$ implies
\begin{equation*}
    \dim(Y_n) \le n + \sup\big\{\hdim^{G}\bigl(X(c)\bigr) \,\big|\, c\in\ob(\calc)\big\}.
\end{equation*}
Having shown this, we beg the reader's pardon for omitting the actual
construction of the $G$-homotopy equivalence $f_n\colon Z_n \times_{\calc} X
\to Y_n$. This construction is based on various standard cofibration arguments
and is essentially the same as in the proof that the $G$-homotopy type of a
$G$-$CW$-complex depends only on the $G$-homotopy classes of the attaching
maps.
\\[1mm]\ref{lem:Z_times_calc_X:fixed_point_sets}
There is an obvious natural map
\begin{equation*}
    f_Z \colon Z \times_{\calc} X^H \to \left(Z \times_{\calc} X\right)^H,
\end{equation*}
and analogously for $Z_n$ in place of $Z$. We show by induction over $n$ that
the map $f_{Z_n}$ is a homeomorphism. Then the claim follows for all
finite-dimensional $\calc$-$CW$-complexes $Z$. The infinite-dimensional case
follows using a colimit argument.

In the induction step, just note that the various homeomorphisms $f_{Z_{n-1}}$
induce a map from the pushout
\begin{equation*}
    \xymatrix{\coprod_{j \in J}X(c_j)^H\times S^{n-1}
        \ar[d]_{\coprod_{j\in J}\id_{X(c_j)^H}\times i} \ar[r]
        & Z_{n-1}\times_\calc X^H \ar[d] \\
        \coprod_{j \in J}X(c_j)^H\times D^n \ar[r] & Z_n\times_\calc X^H
    }
\end{equation*}
to the pushout
\begin{equation*}
    \xymatrix{\Bigl(\coprod_{j \in J} X(c_j) \times S^{n-1}\Bigr)^H
        \ar[d]_{\bigl(\coprod_{j\in J}\id_{X(c_j)}\times i\bigr)^H} \ar[r]
        & (Z_{n-1}\times_\calc X)^H \ar[d] \\
        \Bigl(\coprod_{j \in J}X(c_j)\times D^n\Bigr)^H \ar[r] & (Z_n\times_\calc X)^H
    }
\end{equation*}
and that a pushout of homeomorphisms is again a homeomorphism.
\end{proof}

\begin{lemma} \label{lemma:hocolim_G/im(psi_i)_simeq_pt}
Let $\{G_i \mid i \in I\}$ be a directed system of groups.  Let $G = \colim_{i
\in I} G_i$ be its colimit.  Let $\calf$ be a family of subgroups of $G$ such
that every $H \in \calf$ is contained in $\im(\psi_i)$ for some $i \in I$.  We
obtain a covariant functor $G/\im(\psi_?) \colon I \to G\text{-}\SPACES$ by
sending $i$ to $G/\im(\psi_i)$.

Then $EI \times_I G/\im(\psi_?)$ is $G$-homotopy equivalent to $\pt$.
\end{lemma}
\begin{proof}
We already know from
Lemma~\ref{lem:Z_times_calc_X}~\ref{lem:Z_times_calc_X:G-CW-Structure} that $EI
\times_I G/\im(\psi_?)$ has the $G$-homotopy type of a $G$-$CW$-complex (one
can actually check easily that it is itself a $G$-$CW$-complex, the
$n$-skeleton being $EI_n \times_I G/\im(\psi_?)$). Hence, by the equivariant
Whitehead theorem (see for instance~\cite[Theorem~2.4 on
page~36]{Lueck(1989)}), it suffices to show that $\left(EI \times_I
G/\im(\psi_?)\right)^H$ is weakly contractible for every subgroup $H \subseteq
G$. Because of
Lemma~\ref{lem:Z_times_calc_X}~\ref{lem:Z_times_calc_X:fixed_point_sets}, it
remains to show for every $H\subseteq G$ that $EI \times_I (G/\im(\psi_?))^H$
is weakly contractible.

Recall that $EI \times_I (G/\im(\psi_?))^H$ is defined as $\coprod_{i \in I}
EI(i) \times (G/\im(\psi_i))^H/\sim$ for some tensor product-like equivalence
relation $\sim$ (see~\cite[Definition~1.4]{Davis-Lueck(1998)}).  Consider two
points $x_1$ and $x_2$ in $EI \times_I (G/\im(\psi_?))^H$. Since $G =
\bigcup_{i \in I} \im(\psi_i)$ and, by assumption, $I$ is directed and there is
an index $j \in I$ with $H \subseteq \im(\psi_j)$, we can find an index $k \in
I$ and $z_1, z_2 \in EI(k)$ such that $(z_l,1\cdot \im(\psi_k))\in EI(k) \times
G/\im(\psi_k)$ lies in $EI(k) \times (G/\im(\psi_k))^H$ and represent $x_l$ for
$l = 1,2$.  As $EI(k)$ is contractible and in particular path-connected, we can
join $z_1$ and $z_2$ and hence $x_1$ and $x_2$ by a path. This shows that $EI
\times_I (G/\im(\psi_?))^H$ is path-connected.

It remains to show for every compact subspace $C \subseteq EI \times_I
(G/\im(\psi_?))^H$ that we can find a subspace $D \subseteq EI \times_I
(G/\im(\psi_?))^H$ such that $C\subseteq D$ and such that $D$ is homotopy
equivalent to a discrete set. In the following, we use the functorial bar model
as described in \cite[page~230]{Davis-Lueck(1998)} as the model for $EI$.
Notice that the $n$-cells in $EI$ are indexed by the set of all sequences $i_0
< i_1 < \ldots < i_n$ of elements in $I$. Since $C$ is compact, it is contained
in the union of finitely many of the equivariant cells of $EI \times_I
G/\im(\psi_?)$. Such an equivariant cell is of the shape $(\mor_I(?,i) \times
D^n) \times_I G/\im(\psi_?) = G/\im(\psi_i) \times D^n$ for some cell
$\mor_I(?,i) \times D^n$ of $EI$. Since $I$ is directed, we can find an index
$j \in I$ such that $C$ is already contained in $(EI)_{\le j} \times_I
G/\im(\psi_?)$, where $(EI)_{\le j}$ is the $I$-$CW$-subcomplex of $EI$ built
of cells for which the indexing sequence $i_0 < i_1 < \ldots < i_n$ satisfies
$i_n \le j$. Let $k_j \colon I_{\le j} \to I$ be the inclusion of the full
subcategory $I_{\le j}$ of $I$ formed by the objects $i \in I$ with $i \le j$.
Then $(EI)_{\le j}$ is isomorphic to the induction $(k_j)_*E(I_{\le j})$ of the
bar model $E(I_{\le j})$ of the classifying space for the category $I_{\le j}$.
However, $I_{\le j}$ has a terminal object, namely $j$, so the canonical
projection $E(I_{\le j}) \to \mor_{I_{\le j}}(?,j)$ is a homotopy equivalence
of contravariant $I_{\le j}$-spaces. It induces a homotopy equivalence
\begin{multline*}
  (EI)_{\le j} \times_I G/\im(\psi_?) \cong (k_j)_*E(I_{\le j}) \times_I G/\im(\psi_?)
  \\
  \xrightarrow{\simeq} (k_j)_*\mor_{I_{\le j}}(?,j) \times_I G/\im(\psi_?) =
  G/\im(\psi_j).
\end{multline*}
Summarizing, $C \subseteq D$ holds for $D:= (EI)_{\le j} \times_I
G/\im(\psi_?)$, and $D$ is homotopy equivalent to a discrete set.
\end{proof}

Now we can give the proof of Theorem~\ref{the:model_for_colimits}.
\begin{proof}
Fix a model $\EGF{G}{\calf}$. We get a covariant functor $G/\im(\psi_?) \times
\EGF{G}{\calf} \colon I \to G\text{-}\SPACES$ by sending $i$ to $G/\im(\psi_i)
\times \EGF{G}{\calf}$ and, for $i \le j$, by sending the morphism $i \to j$ in
$I$ to the $G$-map $\pr \times\id \colon G/\im(\psi_i) \times \EGF{G}{\calf}
\to G/\im(\psi_j) \times \EGF{G}{\calf}$, where $\pr$ is the obvious
projection. From Lemma~\ref{lemma:hocolim_G/im(psi_i)_simeq_pt}, we obtain a
$G$-homotopy equivalence
\begin{multline*}
  EI\times_I \left(G/\im(\psi_?) \times \EGF{G}{\calf} \right) \cong \left( EI \times_I
    G/\im(\psi_?) \right)\times \EGF{G}{\calf}
  \\
  \xrightarrow{\simeq} \pt \times \EGF{G}{\calf} = \EGF{G}{\calf}.
\end{multline*}
The $G$-space $G/\im(\psi_i) \times \EGF{G}{\calf}$ is a $G$-$CW$-complex for
every $i$ and is $G$-homeomorphic to $G \times_{G_i} \psi_i^*\EGF{G}{\calf}$,
where $\psi_i^*\EGF{G}{\calf}$ is the $G_i$-$CW$-complex obtained from the
$G$-$CW$-complex $\EGF{G}{\calf}$ by restriction with $\psi_i$. One easily
checks that $\psi_i^*\EGF{G}{\calf}$ is a model for $\EGF{G_i}{\psi_i^*\calf}$.
Hence, the $G$-space $G/\im(\psi_i) \times \EGF{G}{\calf}$ is $G$-homotopy
equivalent to a $G$-$CW$-complex of dimension
$\hdim^{G_i}(\EGF{G_i}{\psi_i^*\calf})$. Now the claim follows from
Lemma~\ref{lem:Z_times_calc_X}~\ref{lem:Z_times_calc_X:G-CW-Structure} applied
to $Z = EI$ and $X = G/\im(\psi_?) \times \EGF{G}{\calf}$.
\end{proof}


\typeout{-------------------- Section 5 --------------------------}

\section{On the dimension of $\edub{G}$}
\label{sec:On_the_dimension_of_edub(G)}

In this section, we deal with the question what can be said about the dimension
of $G$-$CW$-complexes which are models for $\edub{G}$. The same question for
$\eub{G}$ has already been thoroughly investigated, and in many cases
satisfying answers are known (see for instance~\cite{Lueck(2005s)}, where
further references to relevant papers are given).


\subsection{Lower bounds for the homotopy dimension for $\calvcyc$ by the one for $\calfin$}
\label{subsec:Lower_bounds_for_the_homotopy_dimension_for_calvcyc_by_the_one_for_calfin}

In this subsection, we give lower bounds for $\hdim^{G}(\edub{G})$ in terms of
$\hdim^{G}(\eub{G})$. We have introduced the notion of the homotopy dimension
$\hdim^G(X)$ of a $G$-$CW$-complex $X$ in
Definition~\ref{def:homotopy_dimension}.

\begin{proposition} \label{prop:_comparing_calf_subseteq_calg}
Let $\calf\subseteq\calg$ be families of subgroups of $G$.
\begin{enumerate}

\item \label{prop:_comparing_calf_subseteq_calg:mindim}
Let $n \ge 0$ be an integer. Suppose that for any $H\in\calg$ there is an
$n$-dimensional model for $\EGF{H}{\calf\cap H}$. Then
\begin{equation*}
    \hdim^G(\EGF{G}{\calf}) \le n + \hdim^G(\EGF{G}{\calg});
\end{equation*}

\item \label{prop:_comparing_calf_subseteq_calg:finite_type}
Suppose that there exists a finite $H$-$CW$-model for $\EGF{H}{\calf\cap H}$
for every $H\in\calg$ and a finite $G$-$CW$-model for $\EGF{G}{\calg}$. Then
there exists a finite $G$-$CW$-model for $\EGF{G}{\calf}$.

The same is true if we replace ``finite'' by ``of finite type'' everywhere.
\end{enumerate}
\end{proposition}
\begin{proof}
We only give the proof of
assertion~\ref{prop:_comparing_calf_subseteq_calg:mindim} since the one for
assertion~\ref{prop:_comparing_calf_subseteq_calg:finite_type} is similar.

Let $Z$ be an $m$-dimensional $G$-$CW$-complex with isotropy groups in $\calg$.
We will show that then $Z\times E_\calf(G)$ is $G$-homotopy equivalent to an
$(n+m)$-dimensional $G$-$CW$-complex, which implies the claim of the
proposition as $E_\calg(G)\times E_\calf(G)$ is a model for $E_\calf(G)$.

We utilize induction over the dimension $d$ of $Z$. If $Z=\emptyset$, then
there is nothing to show, so let $d\geq 0$. Crossing the $G$-pushout telling
how $Z_d$ arises from $Z_{d-1}$ with $E_\calf(G)$ yields a $G$-pushout
\begin{equation} \label{eq:E_F/G(G)_1}
    \raisebox{2.5em}{\xymatrix{\coprod_{i\in I_d} G/H_i\times \EGF{G}{\calf} \times S^{d-1} \ar[d] \ar[r]^-q &
            Z_{d-1}\times \EGF{G}{\calf}  \ar[d] \\
            \coprod_{i\in I_d}G/H_i\times \EGF{G}{\calf}\times D^d \ar[r] &
            Z_d\times \EGF{G}{\calf}
    }}
\end{equation}
There is a $G$-homotopy equivalence $f_i\colon
G\times_{H_i}E_{\calf\cap H_i}(H_i)\to G/H_i\times E_\calf(G)$
since the restriction of $E_\calf(G)$ to $H_i$ is a $H_i$-$CW$-model for
$E_{\calf\cap H_i}(H_i)$. We put $f_S := \coprod_i
f_i\times\id_{S^{d-1}}$ and $f_D := \coprod_i f_i\times\id_{D^d}$.  Furthermore, by
induction hypothesis, there is a $G$-homotopy equivalence $h\colon Z'\to Z_{d-1}\times
E_\calf(G)$, where $Z'$ is an $(n+d-1)$-dimensional $G$-$CW$-complex. Choose a
$G$-homotopy inverse $k$ of $h$.

For $j\colon\coprod_iG\times_{H_i}\EGF{H_i}{\calf\cap H_i} \times
S^{d-1}\to\cyl(k\circ q\circ f_S)$ and $p\colon\cyl(k\circ q\circ f_S)\to Z'$
the obvious inclusion and projection, $h\circ p\circ j$ and $q\circ f_S$ are
$G$-homotopy equivalent. Since $j$ is a $G$-cofibration, $h\circ p$ can be
altered within its $G$-homotopy class to yield a $G$-map $h'\colon\cyl(k\circ
q\circ f_S)\to Z_{d-1}\times E_\calf(G)$ such that $h'\circ j=q\circ f_S$. Now
consider the $(n+d)$-dimensional $G$-$CW$-complex $Z''$ which is defined by the
$G$-pushout
\begin{equation} \label{eq:E_F/G(G)_2}
    \raisebox{2.5em}{\xymatrix{\coprod_{i\in I_d}G\times_{H_i} \EGF{H_i}{\calf\cap H_i}\times S^{d-1}
        \ar[d] \ar[r]^-j & \cyl(k\circ q\circ f_S) \ar[d] \\
        \coprod_{i\in I_d}G\times_{H_i} \EGF{H_i}{\calf\cap H_i} \times D^d
        \ar[r] & Z''
    }}
\end{equation}
The $G$-homotopy equivalences $f_S$, $f_D$ and $h'$ induce a map of $G$-pushouts from
\eqref{eq:E_F/G(G)_2} to \eqref{eq:E_F/G(G)_1}, and, as the left vertical arrows in these
diagrams are $G$-cofibrations, $Z_d\times \EGF{G}{\calf}$ is $G$-homotopy equivalent to
$Z''$ by \cite[Lemma~2.13]{Lueck(1989)}.
\end{proof}

Since any virtually cyclic group $V$ admits a finite $1$-dimensional model for
$\eub{V}$, Proposition~\ref{prop:_comparing_calf_subseteq_calg} implies:

\begin{corollary}\label{cor:dim(eub(G)_le_dim(edub(G)_plus_1}
Let $G$ be a group.

\begin{enumerate}
\item \label{cor:dim(eub(G)_le_dim(edub(G)_plus_1:mindim}
We have $\hdim^G(\eub{G}) \le 1 + \hdim^G(\edub{G})$;

\item \label{cor:dim(eub(G)_le_dim(edub(G)_plus_1:type}
If there is a model for $\edub{G}$ which is finite or of finite type, then
there is a model for $\eub{G}$ which is finite or of finite type respectively.
\end{enumerate}
\end{corollary}

\begin{remark}
The inequality appearing in
Corollary~\ref{cor:dim(eub(G)_le_dim(edub(G)_plus_1}~\ref{cor:dim(eub(G)_le_dim(edub(G)_plus_1:mindim}
is sharp (see Example~\ref{the:Heisenberg_group_extension_by_Z},
Theorem~\ref{the:low_hdim}~\ref{the:low_hdim:G_locally_virtually_cyclic} and
Example~\ref{exa:hdim_for_eub_and_edub_of_Z[1/p]}).
\end{remark}

\begin{corollary} \label{cor:products}
If $G$ and $H$ are groups, then
\begin{equation*}
    \hdim^{G\times H}(\edub{(G \times H)}) \le \hdim^G(\edub{G})+\hdim^H(\edub{H}) + 3.
\end{equation*}
\end{corollary}
\begin{proof}
Let $\calf$ be the family of subgroups of $G \times H$ which are contained in
subgroups of the shape $V \times W$ for virtually cyclic subgroups $V \subseteq
G$ and $W \subseteq H$.  Then $\edub{G} \times \edub{H}$ is a model for $\EGF{G
\times H}{\calf}$. An element $K$ in $\calf$ contains $\IZ^n$ for some $n \in
\{0,1,2\}$ as subgroup of finite index.  This implies $\hdim^K(\edub{K}) \le 3$
(see Example~\ref{exa:mindim(G)_G_virtually_Zn}).  Now apply
Proposition~\ref{prop:_comparing_calf_subseteq_calg} to $G \times H$ for
$\calvcyc \subseteq \calf$.
\end{proof}

\begin{remark}
The inequality appearing in Corollary~\ref{cor:products} is sharp, as can
already be seen in the case $G=H=\IZ$ (see
Example~\ref{exa:mindim(G)_G_virtually_Zn}).

For another example, let $G_{-1}$ be the group appearing in
Example~\ref{the:Heisenberg_group_extension_by_Z}. Then
$\hdim^{G_{-1}}(\edub{G_{-1}}) = 3$. Since the center of $G_{-1} \times G_{-1}$
is isomorphic to $\IZ^2$ and $\vcd(G_{-1} \times G_{-1}) = \vcd(G_{-1}) +
\vcd(G_{-1}) = 8$ by
Lemma~\ref{lem:prop_virtually_poly_Z-gr}~\ref{lem:prop_virtually_poly_Z-gr:exact_sequences_and_vcd},
Theorem~\ref{the:virtually_poly_Z-groups}~\ref{the:virtually_poly_Z-groups:ge_2}
implies $\hdim^{G_ {-1} \times G_{_1}}(\edub{(G_ {-1} \times G_{-1})}) = 9$.
\end{remark}


\subsection{Upper bounds for the homotopy dimension for $\calvcyc$ by the one for $\calfin$}
\label{subsec:Upper_bounds_for_the_homotopy_dimension_for_calvcyc_by_the_one_for_calfin}

The goal of this subsection is to give upper bounds for $\hdim^{G}(\edub{G})$
in terms of $\hdim^{G}(\eub{G})$. Recall Notation~\ref{not:(Mcalg,calf)}.

\begin{theorem} \label{the:mindim_for_G_with_NM_orM_fin_vcyc}

\begin{enumerate}

\item \label{the:mindim_for_G_with_NM_orM_fin_vcyc:M}
Let $G$ be a group satisfying $M_{\calfin\subseteq\calvcyc}$. Suppose that one of the following
conditions is satisfied:
\begin{enumerate}
\item \label{the:mindim_for_G_with_NM_orM_fin_vcyc:M:condV}
For every maximal infinite virtually cyclic subgroup $V \subseteq G$ 
of type $I$ we have
$$\dim(EW_GV) < \infty.$$

\item \label{the:mindim_for_G_with_NM_orM_fin_vcyc:M:condG}
We have: $\hdim^G(\edub{G}) < \infty$.
\end{enumerate}

Then:
\begin{equation*}
  \hdim^G(\edub{G}) \le \hdim^G(\eub{G}) +1;
\end{equation*}

\item \label{the:mindim_for_G_with_NM_orM_fin_vcyc:NM}
Let $G$ be a group satisfying $N\!M_{\calfin\subseteq\calvcyc}$. Then:
\begin{equation*}
    \hdim^G(\edub{G})
    \begin{cases}
    = \hdim^G(\eub{G}) & \text{if } \hdim^G(\eub{G}) \ge 2; \\
    \le 2 & \text{if } \hdim^G(\eub{G}) \le 1.
    \end{cases}
\end{equation*}

\end{enumerate}
\end{theorem}
\begin{proof}%
\ref{the:mindim_for_G_with_NM_orM_fin_vcyc:M} 
We first explain that condition~\ref{the:mindim_for_G_with_NM_orM_fin_vcyc:M:condG} 
implies condition~\ref{the:mindim_for_G_with_NM_orM_fin_vcyc:M:condV}.
Choose a finite-dimensional $G$-$CW$-model for $\edub{G}$. Let $V\subseteq G$ be
a maximal infinite virtually cyclic subgroup. The $G$-action on $\edub{G}$ 
induces a $N_GV$-action by restriction and  hence a $W_GV$-action on $\edub{G}^V$.
Obviously $\edub{G}^V$ is a $W_GV$-$CW$-complex which has finite dimension and
is contractible (after forgetting the
group action). It suffices to show that it is a model for $EW_GV$, i.e., that $W_GV$ acts 
freely on $\edub{G}^V$. For $x \in \edub{G}^V$ the isotropy group $G_x$ under the 
$G$-action is a virtually cyclic subgroup of $G$ which contains $V$. By the maximality
of $V$ we conclude $G_x = V$. Hence the isotropy group of $x$ under the $W_GV$-action is
trivial. 

Hence it suffices to prove assertion~\ref{the:mindim_for_G_with_NM_orM_fin_vcyc:M}
provided that condition~\ref{the:mindim_for_G_with_NM_orM_fin_vcyc:M:condV} is true.

Let $\cyl(f_V)$ be the
mapping cylinder of the $N_GV$-map $f_V\colon\eub{N_GV}\to EW_GV$
appearing in Corollary~\ref{cor:Efin_and_Evcyc_model_assuming_(McalG,calF)}.
Then we obtain from
Corollary~\ref{cor:Efin_and_Evcyc_model_assuming_(McalG,calF)} a cellular
$G$-pushout
\begin{equation*}
    \xymatrix{\coprod_{V\in\calm} G\times_{N_GV}\eub{N_GV}
        \ar[d]^{\coprod_{V\in\calm} G\times_{N_GV} i_V} \ar[r]^-i & \eub{G} \ar[d] \\
        \coprod_{V\in\calm} G\times_{N_GV} \cyl(f_V) \ar[r] & \edub{G}
    }
\end{equation*}
where the $i_V$ are the obvious inclusions. Since $\res_G^{N_GV}\eub{G}$ is an
$N_GV$-$CW$-model for $\eub{N_GV}$, it suffices to show
\begin{equation} \label{the:mindim_for_G_with_NM_orM_fin_vcyc:(1)}
    \hdim^{W_GV}(EW_GV) \le \hdim^{N_GV}(\eub{N_GV}).
\end{equation} 
Let $d = \hdim^{W_GV}(EW_GV)$ if $\hdim^{W_GV}(EW_GV)$ is finite,
or let $d$ be any positive integer if $\hdim^{W_GV}(EW_GV)$ is infinite. 
We can choose a left $\IZ[W_GV]$-module $M$ such
that $H^d_{W_GV}(EW_GV;M) \not=0$ (see~\cite[Proposition~2.2 on page
185]{Brown(1982)}). Since $W_GV$ is torsionfree, $V\backslash\eub{N_GV}$ is a
free $W_GV$-$CW$-complex. By the equivariant Whitehead theorem (see
\cite[Theorem~2.4 on page~36]{Lueck(1989)}), the projection $V\backslash
\eub{N_GV} \times EW_GV \to V\backslash \eub{N_GV}$ is a $W_GV$-homotopy
equivalence, where $W_GV$ acts diagonally on the source. Since $EW_GV$ is a
free $W_GV$-$CW$-complex, we obtain a cohomology spectral
sequence converging to $H^{p+q}_{W_GV}(V\backslash \eub{N_GV};M)$ whose
$E_2$-term is given by
\begin{equation*}
    E^{p,q}_2 = H^p_{W_GV}(EW_GV;H^q_{W_GV}(V\backslash \eub{N_GV} \times W_GV;M)).
\end{equation*}
Here, $W_GV$ acts diagonally from the left on $V\backslash \eub{N_GV} \times
W_GV$, and the $W_GV$-action on $H^q_{W_GV}(V\backslash \eub{N_GV} \times
W_GV;M))$ comes from the right $W_GV$-action on $V\backslash \eub{N_GV} \times
W_GV$ which is given by $(x,w') \cdot w = (x,w'w)$ (and which commutes with the
left $W_GV$-action). Now consider the homeomorphism
\begin{equation*}
  \phi \colon V\backslash \eub{N_GV} \times W_GV \to V\backslash \eub{N_GV} \times W_GV,
  \quad (x,w) \mapsto (wx, w).
\end{equation*}
It is equivariant with respect to the left $W_GV$-actions given by $w \cdot
(x,w') = (x,ww')$ on the source and the diagonal action on the target.
Moreover, it is equivariant with respect to the right $W_GV$-actions given by
$(x,w') \cdot w = (w^{-1}x,w'w)$ on the source and by $(x,w') \cdot w =
(x,w'w)$ on the target. Thus, we obtain an isomorphism of abelian groups
\begin{multline*}
    \alpha \colon H^q_{W_GV}(V\backslash \eub{N_GV} \times W_GV;M)
    \xrightarrow{H^q_{W_GV}(\phi;M)} H^q_{W_GV}(V\backslash
    \eub{N_GV} \times W_GV;M) \\
    \xrightarrow{\mu} H^q(V\backslash \eub{N_GV};M).
\end{multline*}
where $\mu$ is induced by the isomorphism of $\IZ$-cochain complexes
\begin{multline*}
    \hom_{\IZ[W_GV]}(C_*(V\backslash \eub{N_GV} \times W_GV),M) \\
    \cong\hom_{\IZ[W_GV]}(C_*(V\backslash \eub{N_GV}) \otimes_{\IZ} \IZ[W_GV],M) \cong
    \hom_{\IZ}(C_*(V\backslash \eub{N_GV});M).
\end{multline*}
The isomorphism $\alpha$ becomes an isomorphism of left $\IZ[W_GV]$-modules if
we equip the target with the $W_GV$-action that comes from the following
$W_GV$-action on the $\IZ$-chain complex $\hom_{\IZ}(C_*(V\backslash
\eub{N_GV});M)$: namely, $w \in W_GV$ acts by sending $f$ to $l_w \circ f \circ
C_*(l_{w^{-1}})$, where $l_w$ and $l_{w^{-1}}$ denotes left multiplication with
$w$ and $w^{-1}$. The universal coefficient theorem for the ring $\IZ$ yields
an exact sequence of $\IZ[W_GV]$-modules
\begin{multline*}
    0 \to \operatorname{Ext}^1_{\IZ}(H_{q-1}(V\backslash \eub{N_GV}),M) \to
    H^q(V\backslash \eub{N_GV};M) \\
    \to \hom_{\IZ}(H_q(V\backslash \eub{N_GV}),M) \to 0,
\end{multline*}
where the $W_GV$-actions on the first and third term are given by combining the
obvious right $W_GV$-action on $H_*(V\backslash \eub{N_GV})$ and the given left
$W_GV$-action on $M$.

Now we separately treat the two cases that $V$ is a virtually cyclic group of
type I or of type {II}.  In the first case, we can find an epimorphism $\alpha
\colon V \to \IZ$ with finite kernel.  Then $E\IZ$, viewed as a $V$-space by
restricting with $\alpha$, is a model for $\eub{V}$. Similarly,
$\res_{N_GV}^V\eub N_GV$ is a model for $\eub{V}$. Hence, there is a
$V$-homotopy equivalence $\eub{N_GV} \to E\IZ$, which induces a homotopy
equivalence $V\backslash \eub{N_GV} \to \IZ\backslash E\IZ$.  We conclude that
$H_q(V\backslash \eub{N_GV})$ is infinite cyclic for $q = 0,1$ and trivial for
$q \ge 2$.  This implies that we obtain an isomorphism of $\IZ[W_GV]$-modules
\begin{equation*}
    H^q(V\backslash \eub{N_GV};M) \xrightarrow{\cong}
    \hom_{\IZ}(H_q(V\backslash \eub{N_GV}),M).
\end{equation*}
Let $\epsilon \colon W_GV \to \{\pm 1\}$ the homomorphism which is uniquely
determined by the property that $w \in W_GV$ acts on $H_1(V\backslash
\eub{N_GV})$ by multiplication with $\epsilon(w)$.  Let $M^{\epsilon}$ be the
$\IZ[W_GV]$-module whose underlying abelian group is $M$ but on which the
$W_GV$-action is twisted by $\epsilon$, that is, $w \cdot x$ is defined to be
$\epsilon(w)wx$ for $w \in W_GV$ and $x \in M$. Then
$\hom_{\IZ}(H_1(V\backslash \eub{N_GV}),M)$ is $\IZ[W_GV]$-isomorphic to
$M^{\epsilon}$. Summarizing, what we have shown so far is that $E_2^{p,q} =
H^p_{W_GV}(EW_GV;M^{\epsilon})$ for $q = 1$ and $E_2^{p,q}=0$ for $q \ge 2$.
We get $d < \infty$ from condition~\ref{the:mindim_for_G_with_NM_orM_fin_vcyc:M:condV}.
Obviously $E_2^{p,q} = 0$ holds for $p > d$.
So the $E^2$-term of the first quadrant spectral sequence looks like
$$\xymatrix{
\vdots & \vdots & \vdots & \vdots & \vdots & \vdots & 
\\
0 &   \cdots & 0 & 0 & 0 & 0 &\cdots
\\
0 &  \cdots & 0 & 0 & 0 & 0 &\cdots
\\
E^{0,1}_2   & \cdots & E^{2,d-1}_2 
& H^d_{W_GV}(EW_GV;M^{\epsilon}) & 0 & 0 &\cdots
\\
E^{0,0}_2    & \cdots & E^{d-1,0}_2& E^{d,0}_2 & 0 & 0 & \cdots
}$$
Hence we obtain an isomorphism
\begin{equation*}
    H^{d+1}_{W_GV}(V\backslash \eub{N_GV};M) \cong H^d_{W_GV}(EW_GV;M^{\epsilon}).
\end{equation*}
Since $(M^{\epsilon})^{\epsilon} = M$ and $H^d_{W_GV}(EW_GV;M) \not= 0$, this
implies that we have $H^{d+1}_{W_GV}(V\backslash \eub{N_GV};M^{\epsilon}) \not=
0$. It follows that $\hdim^{N_GV}(\eub{N_GV}) \ge d+1$, and hence
assertion~\eqref{the:mindim_for_G_with_NM_orM_fin_vcyc:(1)} is true in the case
that $V$ is of type I.

It remains to treat the (easier) case {II}, in which there is an epimorphism
$\alpha \colon V \to D_{\infty}$ onto the infinite dihedral group. Just as
above, we obtain a homotopy equivalence $V\backslash \eub{N_GV} \to
D_{\infty}\backslash\eub{D_{\infty}}$, and as
$D_{\infty}\backslash\eub{D_{\infty}}$ is contractible, $H_q(V\backslash
\eub{N_GV})$ is trivial for $q \ge 1$ and is isomorphic as a $\IZ[W_GV]$-module
to $\IZ$ with the trivial $W_GV$-action for $q = 0$. Hence, the cohomology
spectral sequence yields an isomorphism 
\begin{equation*}
    H^d_{W_GV}(V\backslash \eub{N_GV};M) \cong H^d_{W_GV}(EW_GV;M),
\end{equation*}
whence it follows that $H^d_{W_GV}(V\backslash \eub{N_GV};M) \not= 0$. So
$\hdim^{N_GV}(\eub{N_GV}) \ge d$, and we get
assertion~\eqref{the:mindim_for_G_with_NM_orM_fin_vcyc:(1)}.  
\\[1mm]\ref{the:mindim_for_G_with_NM_orM_fin_vcyc:NM}
From Corollary~\ref{cor:Efin_and_Evcyc_model_with_N_GM_is_M} we get a cellular
$G$-pushout
\begin{equation*}
    \xymatrix{\coprod_{V\in\calm} G\times_{V}\eub{V}
        \ar[d]^{\coprod_{V\in\calm} G\times_{V} i_V} \ar[r]^-i & \eub{G} \ar[d] \\
        \coprod_{V\in\calm} G\times_{V} \cone(\eub{V}) \ar[r] & \edub{G}
}
\end{equation*}
where $\calm$ is a complete system of representatives of the conjugacy classes
of maximal infinite virtually cyclic subgroups and $i_V$ is the obvious
inclusion. Since each $V$ is virtually cyclic, we can choose a $1$-dimensional
model for $\eub{V}$. Then $\cone(\eub{V})$ is $2$-dimensional. This implies
that there exists a $\max\{d,2\}$-model for $\edub{G}$ if there is a
$d$-dimensional model for $\eub{G}$. Hence it remains to prove
$\hdim^G(\eub{G}) \le \hdim^G(\edub{G})$ provided that $\hdim^G(\eub{G}) \ge 2$
holds.

Put $d = \hdim^G(\eub{G})$. We get from~\cite[Proposition~11.10 on
page~221]{Lueck(1989)} a contravariant $\IZ\OrGF{G}{\calfin}$-module $M$ with
$H^d\bigl(\hom_{\IZ\OrGF{G}{\calfin}}(C_*(\eub{G}),M)\bigr) \not= 0$. Define a
contravariant $\IZ\OrGF{G}{\calvcyc}$-module $\overline{M}$ by putting
$\overline{M}(G/H) = M(G/H)$ if $H \in \calfin$ and $\overline{M}(G/H) = 0$ if
$H \notin \calfin$. If $f \colon G/H \to G/K$ is a morphism in
$\OrGF{G}{\calvcyc}$, define $\overline{M}(f) = M(f)$ if $H,K \in \calfin$, and
$\overline{M}(f) = 0$ otherwise. If $j \colon \OrGF{G}{\calfin} \to
\OrGF{G}{\calvcyc}$ is the inclusion, then the restriction $j^*\overline{M}$
agrees with $M$. The adjunction between induction and restriction yields an
isomorphism of $\IZ$-chain complexes
\begin{multline*}
    \hom_{\IZ\OrGF{G}{\calfin}}(C_*(\eub{G}),M)) \cong
    \hom_{\IZ\OrGF{G}{\calfin}}(C_*(\eub{G}),j^*\overline{M})) \\ \cong
    \hom_{\IZ\OrGF{G}{\calvcyc}}(j_*C_*(\eub{G}),\overline{M})).
\end{multline*}
Since there exists a $1$-dimensional model for $\eub{V}$ if $V \in\calm$ and
since $d \ge 2$, we obtain an epimorphism
\begin{equation*}
    H^d\bigl(\hom_{\IZ\OrGF{G}{\calvcyc}}(C_*(\edub{G}),\overline{M})\bigr) \to
    H^d\bigl(\hom_{\IZ\OrGF{G}{\calvcyc}}(j_*C_*(\eub{G}),\overline{M})\bigr)
\end{equation*}
from the long exact cohomology sequence of the pair $(\edub{G},\eub{G})$ with
coefficients in $\overline{M}$. This implies
$H^d\bigl(\hom_{\IZ\OrGF{G}{\calvcyc}}(C_*(\edub{G}),\overline{M})\bigr) \not=
0$ and hence $d \le \hdim^G(\edub{G})$. This finishes the
proof of Theorem~\ref{the:mindim_for_G_with_NM_orM_fin_vcyc}.
\end{proof}

\begin{remark} Condition~\ref{the:mindim_for_G_with_NM_orM_fin_vcyc:M:condV} appearing
Theorem~\ref{the:mindim_for_G_with_NM_orM_fin_vcyc} could be avoided if one could show the 
following: If $1 \to C \to \Gamma \to Q \to 1$ is an extension of groups
such that $C$ is infinite cyclic, $Q$ is torsionfree and there exists a finite-dimensional
model for $B\Gamma$, then there exists a finite-dimensional model for $BQ$.
However, it is unlikely that such a statement is true, although we have no explicit
counterexample. 

Suppose that we have such an extension for which no finite-dimensional model for
$BQ$ exists. Then the Gysin sequence,
whose existence follows from spectral sequence argument
appearing in Theorem~\ref{the:mindim_for_G_with_NM_orM_fin_vcyc},  implies
that for every $\IZ Q$-module $M$ and $n \ge \dim(B\Gamma)$ 
we have $H^{n}_Q(EQ;M) \cong H^{n+2}_Q(EQ;M)$, i.e., the cohomology of $Q$ is for large
$n$ two-periodic with respect to any $\IZ Q$-coefficient module $M$.
\end{remark}

\begin{remark}
The inequality appearing in
Theorem~\ref{the:mindim_for_G_with_NM_orM_fin_vcyc}~\ref{the:mindim_for_G_with_NM_orM_fin_vcyc:M}
is sharp. For instance, we have $\hdim^{\IZ^n}(\eub{\IZ^n}) = n$ and
$\hdim^{\IZ^n}(\edub{\IZ^n}) = n+1$ for $n \ge 2$ (see
Example~\ref{exa:mindim(G)_G_virtually_Zn}).
\end{remark}

\begin{example}[Counterexample to a conjecture of Connolly-Fehr\-mann-Hart\-glass]
\label{exa:Counterexample_to_a_conjecture_of_Connolly-Fehrmann-Hartglass}
Let $G$ be the fundamental group of a closed hyperbolic manifold $M$ of
dimension $d$. Then $G$ is torsionfree and satisfies
$N\!M_{\calfin\subseteq\calvcyc}$ by Example~\ref{exa:word-hyp}. This implies,
using
Theorem~\ref{the:mindim_for_G_with_NM_orM_fin_vcyc}~\ref{the:mindim_for_G_with_NM_orM_fin_vcyc:NM}
\begin{equation*}
    d = \hdim^{\{1\}}(M) = \hdim^G(EG) = \hdim^G(\edub{G}).
\end{equation*}
Obviously, $G$ is a discrete cocompact subgroup of the isometry group of the
$d$-dimensional hyperbolic space, which is a Lie group with only finitely many
path components.

Thus, we get counterexamples to a conjecture stated by
Connolly-Fehrmann-Hartglass~\cite[Conjecture~8.1]{Connolly-Fehrmann-Hartglass(2006)},
which predicts that $\hdim^G(\edub{G})=d+1$.
\end{example}


\subsection{Virtually poly-$\IZ$ groups}
\label{subsec:Virtually_poly-Z_groups}

In this subsection, we give a complete calculation of $\hdim^G(\edub{G})$ for
virtually poly-$\IZ$ groups.  Surprisingly, it will depend \emph{not} only on
$\vcd(G) = \hdim^G(\eub{G})$.

A group $G$ is \emph{poly-$\IZ$} if there is a finite sequence $\{1\} = G_0
\subseteq G_1 \subseteq \ldots \subseteq G_n = G$ of subgroups such that
$G_{i-1}$ is normal in $G_{i}$ with infinite cyclic quotient $G_{i}/G_{i-1}$
for $i = 1,2, \ldots , n$.

We have introduced notations like $[\calvcyc\setminus \calfin]$, $[V]$ and
$N_G[V]$ at the beginning of Section~\ref{sec:Passing_to_larger_families}. In
particular, recall Definition~\ref{def:for_sim}.

\begin{theorem}\label{the:virtually_poly_Z-groups}
Let $G$ be a virtually poly $\IZ$ group which is not virtually cyclic.  Let
$[\calvcyc \setminus \calfin]_f$ be the subset of $[\calvcyc\setminus \calfin]$
consisting of those elements $[V]$ for which $N_G[V]$ has finite index in $G$.
Let $[\calvcyc \setminus \calfin]_f/G$ be its quotient under the $G$-action
coming from conjugation. Then
\begin{equation*}
    \hdim^G(\eub{G}) = \vcd(G),
\end{equation*}
and precisely one of the following cases occurs:

\begin{enumerate}

\item \label{the:virtually_poly_Z-groups:empty}
The set $[\calvcyc\setminus \calfin]_f/G$ is empty. This is equivalent to the
condition that every finite index subgroup of $G$ has a finite center. In this
case,
\begin{equation*}
    \hdim^G(\edub{G}) = \vcd(G);
\end{equation*}

\item \label{the:virtually_poly_Z-groups:1_and_N_GV_is_G}
The set $[\calvcyc \setminus \calfin]_f/G$ contains exactly one element
$[V]\cdot G$.  This is equivalent to the condition that there exists an
infinite normal cyclic subgroup $C \subseteq G$ and for every infinite cyclic
subgroup $D \subseteq G$ with $[G:N_GD] < \infty$ we have $C \cap D \not=
\{1\}$. Then there is the following dichotomy:
\begin{enumerate}
\item \label{the:virtually_poly_Z-groups:1_and_N_GV_is_G:vcd-1}
For every $[W] \cdot G \in [\calvcyc \setminus \calfin]/G$ such that $[W] \cdot
G \not= [V] \cdot G$ we have $\vcd(N_G[W]) \le \vcd(G) -2$. In this case,
$\vcd(G) \ge 4$ and
\begin{equation*}
    \hdim^G(\edub{G}) = \vcd(G) - 1;
\end{equation*}

\item \label{the:virtually_poly_Z-groups:1_and_N_GV_is_G:vcd}
There exists $[W] \cdot G \in [\calvcyc \setminus \calfin]/G$ such that $[W]
\cdot G \not= [V] \cdot G$ and $\vcd(N_G[W]) = \vcd(G) -1 $. In this case,
$\vcd(G) \ge 3$ and
\begin{equation*}
    \hdim^G(\edub{G}) = \vcd(G);
\end{equation*}

\end{enumerate}

\item \label{the:virtually_poly_Z-groups:ge_2}
The set $[\calvcyc \setminus \calfin]_f/G$ contains more than one element. This
is equivalent to the condition that there is a finite index subgroup of $G$
whose center contains a subgroup isomorphic to $\IZ^2$. In this case,
\begin{equation*}
    \hdim^G(\edub{G}) = \vcd(G) + 1.
\end{equation*}

\end{enumerate}
\end{theorem}

The proof of Theorem~\ref{the:virtually_poly_Z-groups} needs some preparations.

\begin{lemma} \label{lem:prop_virtually_poly_Z-gr}

\begin{enumerate}

\item \label{lem:prop_virtually_poly_Z-gr:inheritance}
Subgroups and quotients of poly-$\IZ$ groups are again poly-$\IZ$.  If $1 \to
G_0 \to G_1 \to G_2 \to 1$ is an extension of groups and $G_0$ and $G_2$ are
poly-$\IZ$ groups, then $G$ is a poly-$\IZ$ group.

The same statements are true if one replaces ``poly-$\IZ$'' by ``virtually
poly-$\IZ$'' everywhere;

\item \label{lem:prop_virtually_poly_Z-gr:top_homology}
Let $G$ be a poly-$\IZ$ group of cohomological dimension $\cd(G) = d$. Then
$H^d(G)$ is isomorphic to $\IZ$ or to $\IZ/2$;

\item \label{lem:prop_virtually_poly_Z-gr:vcd_and_Hirsch_rank}
Let $G$ be a poly-$\IZ$ group. Suppose that there is a finite sequence $\{1\} =
G_0 \subseteq G_1 \subseteq \ldots \subseteq G_n = G$ of subgroups such that
$G_{i-1}$ is normal in $G_{i}$ with infinite cyclic quotient $G_{i}/G_{i-1}$
for $i = 1,2, \ldots , n$. Then
\begin{equation*}
    n = \cd(G);
\end{equation*}

\item \label{lem:prop_virtually_poly_Z-gr:exact_sequences_and_vcd}
If $1 \to G_0 \to G_1 \to G_2 \to 1$ is an extension of virtually poly-$\IZ$
groups, then
\begin{equation*}
    \vcd(G_1) = \vcd(G_0) + \vcd(G_2);
\end{equation*}

\item \label{lem:prop_virtually_poly_Z-gr:property_(A)}
A virtually poly-$\IZ$ group $G$ has the following property:
      \begin{quote} $(A)$ \quad \begin {minipage}{90mm} If $\phi \colon G
          \to G$ is an automorphism, $g \in G$ is an element of infinite order and $a,b
          \in \IZ$ are integers such that $\phi(g^a) = g^b$, then $a = \pm b$.\end{minipage}
      \end{quote}

\end{enumerate}
\end{lemma}
\begin{proof}\ref{lem:prop_virtually_poly_Z-gr:inheritance} The elementary proof of this
assertion is left to the reader.
\\[1mm]\ref{lem:prop_virtually_poly_Z-gr:top_homology} We use induction over the number
$n$ for which there is a finite sequence $\{1\} = G_0 \subseteq G_1 \subseteq
\ldots \subseteq G_n = G$ of subgroups such that $G_{i-1}$ is normal in $G_{i}$
with infinite cyclic quotient $G_{i}/G_{i-1}$ for $i = 1,2, \ldots , n$.  If $n
= 0$, then $G$ is trivial and the claim is obviously true. The induction step
from $n-1$ to $n \ge 1$ is done as follows.

By the induction hypothesis, the claim is true for $G_{n-1}$. This means that
$\cd(G_{n-1}) = n-1$, that $H^{n-1}(G_{n-1})$ is isomorphic to $\IZ$ or $\IZ/2$
and that $H^i(G_{n-1}) = 0$ for $i \ge n$. We have the extension $1 \to G_{n-1}
\to G \to \IZ \to 1$ if we identify $G/G_{n-1}$ with $\IZ$. Let $f\colon
G_{n-1}\to G_n$ be the group automorphism induced by conjugation with some
preimage in $G$ of a generator of $\IZ$. The associated Hochschild-Serre
cohomology spectral sequence yields the exact Wang sequence
\begin{equation*}
    H^{n-1}(G_{n-1}) \xrightarrow{\id-H^{n-1}(f)} H^{n-1}(G_{n-1}) \to H^n(G)
    \to H^n(G_{n-1})=0.
\end{equation*}
Since $H^{n-1}(G_{n-1})$ is isomorphic to $\IZ$ or $\IZ/2$, the automorphism
$H^{n-1}(f)$ is $\pm \id$.  Hence $H^n(G) \cong \IZ$ if and only if
$H^{n-1}(G_{n-1}) \cong \IZ$ and $H^{n-1}(f) = \id$, and $H^n(G) \cong \IZ/2$
if and only if either $H^{n-1}(G_{n-1}) \cong \IZ$ and $H^{n-1}(f) = - \id$ or
$H^{n-1}(G_{n-1}) \cong \IZ/2$.
\\[1mm]\ref{lem:prop_virtually_poly_Z-gr:vcd_and_Hirsch_rank}
We use induction over the number $n$ as above.  Since $H^n(G)$ is non-trivial,
we have $\cd(G) \ge n$. On the other hand, the cohomological dimension is
subadditive under extensions, i.e., $\cd(G) \le \cd(G_{n-1}) + \cd(\IZ) = (n-1)
+ 1 = n$. We conclude $\cd(G) = n$.
\\[1mm]\ref{lem:prop_virtually_poly_Z-gr:exact_sequences_and_vcd}
One easily constructs an exact sequence of poly-$\IZ$ groups $1 \to G_0' \to
G_1' \to G_2' \to 1$ such that $G_i'$ is a subgroup of finite index in $G_i$
for $i = 0,1,2$. Now, it suffices to show that $\cd(G_1') = \cd(G_0') +
\cd(G_2')$. By induction over the number $n$ as above but with $G_2'$ in place
of $G$, the claim reduces to the case that $G_2' = \IZ$, which follows from
assertion~\ref{lem:prop_virtually_poly_Z-gr:vcd_and_Hirsch_rank}.
\\[1mm]\ref{lem:prop_virtually_poly_Z-gr:property_(A)}
Consider a group extension $1 \to G_0 \xrightarrow{i} G_1 \xrightarrow{p} G_2
\to 1$ such that $G_0$ is characteristic in $G_1$ and $G_0$ and $G_2$ have
property $(A)$. Then $G_1$ has property $(A)$ by the following argument.
Consider an automorphism $\phi_1 \colon G_1 \to G_1$. It induces automorphisms
$\phi_0 \colon G_0 \to G_0$ and $\phi_2 \colon G_2 \to G_2$ satisfying $\phi_1
\circ i = i \circ \phi_0$ and $\phi_2 \circ p = p \circ \phi_1$ since $G_0$ is
characteristic by assumption. Let $g \in G$ be of infinite order and $a,b \in
\IZ$ such that $\phi_1(g^a) = g^b$. Then $\phi_2(p(g)^a) = p(g)^b$, and if
$p(g)$ has infinite order, we see that $a = \pm b$ since $G_2$ has property
$(A)$ by assumption. If $p(g)$ has finite order, we can choose $c \in \IZ, c
\not= 0 $ such that $g^c \in G_0$. Then $\phi_1(g^a) = g^b$ implies
$\phi_0((g^c)^a) = (g^{c})^b$ in $G_0$ and hence $a = \pm b$ since $G_0$ has
property $(A)$ by assumption. Hence, $G_1$ has property $(A)$.

Now we prove by induction over $\vcd(G)$ that a virtually poly-$\IZ$ group $G$
has property $(A)$. The induction beginning $\vcd(G) = 0$ is trivial since then
$G$ is finite. The induction step is done as follows.  Let $K \subseteq G$ be a
subgroup of finite index such that $K$ is poly-$\IZ$. Then $K_0 = \bigcap
\{\phi(K) \mid \phi \in \aut(G)\}$ is a characteristic subgroup of $G$. Since
$G$ is finitely generated and, thus, contains only finitely many subgroups of
any given index, the intersection defining $K_0$ is an intersection of finitely
many groups of finite index. Hence, $K_0$ is a characteristic subgroup of
finite index in $G$. Since every finite group has property $(A)$, the group $G$
has property $(A)$ if $K_0$ has. Consider the group extension $1 \to K_1 \to
K_0 \to H_1(K_0)/\tors(H_1(K_0)) \to 1$. Obviously, $K_1$ is a characteristic
subgroup of $K_0$. Since $K_0$, being a subgroup of the poly-$\IZ$ group $K$,
is poly-$\IZ$, the abelian group $H_1(K_0)/\tors(H_1(K_0)$ is isomorphic to
$\IZ^n$ for some $n \ge 1$. This implies that $K_1$ is a poly-$\IZ$ group
satisfying $\vcd(K_1) < \vcd(K_0) = \vcd(G)$. By induction hypothesis, $K_1$
has property $(A)$. Clearly, $\IZ^n$ also has property $(A)$. This implies that
$K_0$ and hence $G$ have property $(A)$. This finishes the proof of
Lemma~\ref{lem:prop_virtually_poly_Z-gr}.
\end{proof}

The statements about the top homology of a poly-$\IZ$ group $G$ appearing in
Lemma~\ref{lem:prop_virtually_poly_Z-gr} would be obvious if there were closed
manifold models for $BG$. Instead of dealing with this difficult question we
have preferred to give an easy homological argument.

\begin{lemma} \label{lem:N_G[V]_is_N_GV}
Let $G$ be a virtually poly-$\IZ$ group. Let $V \subseteq G$ be an infinite
virtually cyclic group.  Then we can find an infinite cyclic subgroup $C
\subseteq V$ such that $[V] = [C] \in [\calvcyc\setminus\calfin]$ and $N_GC =
N_G[V]$ holds, and such that a model for $\EGF{N_G[V]}{\calvcyc[V]}$ is given
by $\eub{W_GC}$ considered as an $N_G[V] = N_GC$-$CW$-complex by restriction
with the canonical projection $N_GC \to W_GC$.
\end{lemma}
\begin{proof}
We can assume that $V$ itself is infinite cyclic (otherwise, we pass to a
subgroup of finite index). Since $N_G[V]$ is a subgroup of the virtually
poly-$\IZ$ group $G$, it is virtually poly-$\IZ$ and, in particular, finitely
generated.  Let $\{g_1, g_2, \ldots, g_r\}$ be a finite generating set. Fix a
generator $c \in V$.  Since $g_i \in N_G[V]$, we can find integers $a_i, b_i
\in \IZ$ different from zero such that $g_ic^{a_i}g_i^{-1} = c^{b_i}$. Property
$(A)$ (see
Lemma~\ref{lem:prop_virtually_poly_Z-gr}~\ref{lem:prop_virtually_poly_Z-gr:property_(A)})
implies $b_i = \pm a_i$. Put $a = a_1 a_2 \cdots a_r$.  Then $g_ic^{a}g_i^{-1}
= c^{\pm a}$ holds for $i = 1,2, \ldots, r$. Let $C$ be the infinite cyclic
subgroup of $V$ generated by $c^a$.  Then, of course, $[C] = [V]$. Furthermore,
$g_i \in N_GC$ for $i = 1,2, \ldots, r$, which implies $N_G[V] \subseteq N_GC$.
Hence, $N_G[V] = N_GC$.

It remains to prove that $\eub{W_GC}$ is a model for $\EGF{N_GC}{\calvcyc[C]}$.
This follows from the fact that, under the projection $N_GC\to W_GC$, the
family of all virtually cyclic subgroups of $N_GC$ whose intersection with $C$
is infinite coincides with the family of all finite subgroups of $W_GC$.
\end{proof}

Now we are ready to prove Theorem~\ref{the:virtually_poly_Z-groups}.

\begin{proof}
Because of Lemma~\ref{lem:N_G[V]_is_N_GV}, we can choose a set $I$ of infinite
cyclic subgroups of $G$ and models for $\eub{N_GC}$ and $\eub{W_GC}$ for $C \in
I$ such that the following properties hold:

\begin{itemize}

\item For every element $\eta \in [\calvcyc\setminus \calfin]/G$, there exists precisely
  one $C \in I$ such that $[C] \cdot G = \eta$;

\item $N_G[C] = N_GC$ for all $C\in I$;

\item For all $C\in I$, a model for $\EGF{N_G[C]}{\calvcyc[C]}$ is given by $\eub{W_GC}$,
considered as an $N_GC$-$CW$-complex by restriction with the canonical
projection $N_GC \to W_GC$.

\end{itemize}

From Theorem~\ref{the:passing_from_calf_to_calg} we obtain a $G$-pushout
\begin{equation}
    \raisebox{2.5em}{\xymatrix{\coprod_{C \in I} G\times_{N_GC} \eub{N_GC}
        \ar[d]^{\coprod_{C \in I}\id_G \times_{N_GC} f_C} \ar[r]^-i & \eub{G} \ar[d] \\
        \coprod_{C\in I} G\times_{N_GC} \eub{W_GC} \ar[r] & \edub{G}
    }}
\label{basic_G-pushout}
\end{equation}
such that $f_C$ is a cellular $N_GC$-map for every $C \in I$ and $i$ is an
inclusion of $G$-$CW$-complexes.

In the remainder of the proof, we will frequently use that
Lemma~\ref{lem:prop_virtually_poly_Z-gr}~\ref{lem:prop_virtually_poly_Z-gr:exact_sequences_and_vcd}
implies
\begin{align*}
    \vcd(W_GC) &= \vcd(N_GC) - 1; \\
    \vcd(N_GC) &\le \vcd(G); \\
    \vcd(N_GC) = \vcd(G) \,&\!\Leftrightarrow [G:N_GC] < \infty.
\end{align*}
Furthermore, we conclude from~\cite[Example~5.26 on page 305]{Lueck(2005s)}
that
\begin{align*}
  \hdim^G(\eub{G}) &= \vcd(G); \\
  \hdim^G(\eub{N_GC}) &= \vcd(N_GC) &\mspace{-100mu} \text{for } C \in I; \\
  \hdim^G(\eub{W_GC}) &= \vcd(N_GC)-1 &\mspace{-100mu} \text{for } C \in I.
\end{align*}%
\ref{the:virtually_poly_Z-groups:empty} As $[\calvcyc\setminus \calfin]_f/G =
\emptyset$, we must have $\hdim^G(\eub{W_GC}) < \hdim^G(\eub{N_GC}) < \vcd(G)$,
whereas $\hdim^G(\eub{G})=\vcd(G)$. Hence, we conclude from the 
$G$-pushout~\eqref{basic_G-pushout} together with
Remark~\ref{rem:how_to_apply_theorem_the:passing_from_calf_to_calg} that
$\hdim^G(\edub{G}) \le \vcd(G)$. In order to show $\hdim^G(\edub{G}) \ge
\vcd(G)$ we can assume without loss of generality that $G$ is poly-$\IZ$ since
there is a poly-$\IZ$ subgroup $G'$ of finite index and the restriction of
$\edub{G}$ to $G'$ is a model for $\edub{G'}$. The $G$-pushout
\eqref{basic_G-pushout} yields a pushout of $CW$-complexes
\begin{equation} \label{quotient_of_basic_G-pushout}
    \raisebox{2.5em}{\xymatrix{\coprod_{C \in I} N_GC\backslash\eub{N_GC}
        \ar[d]^{\coprod_{C \in I} N_GC\backslash f_C} \ar[r]^-i & G\backslash \eub{G} \ar[d] \\
        \coprod_{C\in I} W_GC\backslash\eub{W_GC} \ar[r] & G\backslash\edub{G}
    }}
\end{equation}
From the associated long exact Mayer-Vietoris sequence, we obtain an
epimorphism $H^{\vcd(G)}(G\backslash\edub{G}) \to
H^{\vcd(G)}(G\backslash\eub{G})$. However, $H^{\vcd(G)}(G\backslash\eub{G})$ is
non-trivial by
Lemma~\ref{lem:prop_virtually_poly_Z-gr}~\ref{lem:prop_virtually_poly_Z-gr:top_homology},
so $H^{\vcd(G)}(G\backslash\edub{G})$ is also non-trivial. This means that
$\hdim^G(\edub{G}) \ge \vcd(G)$.

As $C_GC$ has finite index in $N_GC$ for $C \in I$, the condition
$[\calvcyc\setminus \calfin]_f/G = \emptyset$ is equivalent to the condition
that for every infinite cyclic subgroup $C \subseteq G$ we have $[G:C_GC] =
\infty$.  Obviously, this is true if and only if the center of any finite index
subgroup of $G$ is finite.  This finishes the proof of
assertion~\ref{the:virtually_poly_Z-groups:empty}.
\\[1mm]\ref{the:virtually_poly_Z-groups:ge_2}
We conclude from the $G$-pushout~\eqref{basic_G-pushout} together with
Remark~\ref{rem:how_to_apply_theorem_the:passing_from_calf_to_calg} that
$\hdim^G(\edub{G}) \le \vcd(G)+1$ as $\hdim^G(\eub{N_GC})$,
$\hdim^G(\eub{W_GC})$ and $\hdim^G(\eub{G})$ are less or equal to $\vcd(G)$. It
remains to show that $\hdim^G(\edub{G}) \ge \vcd(G) + 1$, and it suffices to do
this for a subgroup of finite index.  Hence, we can assume without loss of
generality that $G$ is poly-$\IZ$.  We obtain from the
pushout~\eqref{quotient_of_basic_G-pushout} the short exact sequence
\begin{equation*}
    H^{\vcd(G)}(G\backslash \eub{G}) \to \prod_{C\in I} H^{\vcd(G)}(N_GC\backslash
    \eub{N_GC}) \to H^{\vcd(G)+1}(G\backslash \edub{G}) \to 0.
\end{equation*}
If we can show that $H^{\vcd(G)}(G\backslash \eub{G}) \to \prod_{C\in I}
H^{\vcd(G)}(N_GC\backslash \eub{N_GC}) $ is not surjective, then
$H^{\vcd(G)+1}(G\backslash \edub{G})$ is non-trivial, implying $\hdim^G(\edub{G})
\ge \vcd(G) +1$. By
Lemma~\ref{lem:prop_virtually_poly_Z-gr}~\ref{lem:prop_virtually_poly_Z-gr:top_homology},
$H^{\vcd(G)}(G\backslash \eub{G})$ and $H^{\vcd(N_GC)}(N_GC\backslash
\eub{N_GC})$ are non-trivial cyclic groups. Thus, $H^{\vcd(G)}(G\backslash
\eub{G}) \to \prod_{C\in I} H^{\vcd(G)}(N_GC\backslash \eub{N_GC}) $ cannot be
surjective if there is more than one element $C \in I$ with $\vcd(N_GC) =
\vcd(G)$, which is indeed the case here since, by assumption,
$[\calvcyc\setminus\calfin]_f/G$ contains more than one element.

We still want to show that this latter condition is equivalent to the condition
that there is a finite index subgroup of $G$ whose center contains a copy of
$\IZ^2$. Suppose first that $[\calvcyc \setminus \calfin]_f/G$ contains two
different elements $[V_1]\cdot G$ and $[V_2] \cdot G$.
Lemma~\ref{lem:N_G[V]_is_N_GV} shows that we may assume that $V_i$ is infinite
cyclic and $N_GV_i=N_G[V_i]$ for $i=1,2$. Then, as $C_GV_i$ has finite index in
$N_GV_i$ for $i=1,2$, the subgroup $C_GV_1\cap C_GV_2\subseteq G$ is also of
finite index. Choose infinite cyclic groups $C_1$ and $C_2$ such that $C_1
\subseteq V_1 \cap C_GV_2$ and $C_2 \subseteq V_2 \cap C_GV_1$, and let $H
\subseteq G$ be the subgroup generated by $C_1$ and $C_2$. Then $C_GH$ contains
$C_GV_1 \cap C_GV_2$. Hence $C_GH$ has finite index in $G$ and $H \subseteq
C_GH$.  In particular, $H$ is abelian. Since $[V_1] \not= [V_2]$, we have $C_1
\cap C_2 = \{1\}$. So $H$ is isomorphic to $\IZ^2$, and we have found a
finite index subgroup of $G$, namely $C_GH$, whose center contains $\IZ^2$.

Conversely, suppose that $G' \subseteq G$ is a subgroup of finite index whose
center contains $\IZ^2$. By passing to a subgroup of finite index in $G'$, we
can arrange for $G' \subseteq G$ to be torsionfree and normal with finite
quotient $G/G'$. Since $\cent(G')$ is isomorphic to $\IZ^n$ for some $n \ge 2$,
it contains infinitely many maximal infinite cyclic subgroups. As $G/G'$ is
finite, $[\calvcyc \setminus \calfin]_f/G$ must be infinite and, in particular,
contains more than one element.
\\[1mm]\ref{the:virtually_poly_Z-groups:1_and_N_GV_is_G}
We claim that $[\calvcyc\setminus \calfin]_f/G$ consists of precisely one
element if and only if $G$ contains an infinite cyclic normal subgroup $C$ and
for every infinite cyclic subgroup $D \subseteq G$ with $[G:N_GD] < \infty$ we
have $C \cap D \not= \{1\}$. The ``if''-statement is obvious, so let us show
that $G$ contains an infinite cyclic normal subgroup if $[\calvcyc\setminus
\calfin]_f/G$ consists of precisely one element $[V] \cdot G$. Again, just as
in the proof of~\ref{the:virtually_poly_Z-groups:ge_2}, we can assume by
Lemma~\ref{lem:N_G[V]_is_N_GV} that $V$ is infinite cyclic and
$[G:C_GV]<\infty$. This implies that $K := \bigcap_{g \in G} gC_GVg^{-1}$ is a
normal subgroup of $G$ of finite index. Thus, the center $\cent(K)$ of $K$,
being a characteristic subgroup of $K$, is also normal in $G$. Moreover, we
conclude from assertion~\ref{the:virtually_poly_Z-groups:ge_2} that $\cent(K)$
does not contain $\IZ^2$ as subgroup of finite index. Hence, it is isomorphic
to $\IZ \times A$ for some finite abelian subgroup $A$. If $a$ is the order of
$A$, then $C:= \{a \cdot x \mid x \in \cent(K)\}$ is a characteristic infinite
cyclic subgroup of $\cent(K)$ and hence normal in $G$. The claim follows.

In the remainder, let $C_0 \subseteq G$ be an infinite cyclic normal subgroup.
This means that $[\calvcyc\setminus \calfin]_f/G = \{[C_0] \cdot G\}$. Since we
assume that $G$ is not virtually cyclic, $\vcd(G) \ge 2$. In fact, we have
$\vcd(G) > 2$.  Namely, if $\vcd(G) = 2$, then $G$ contains a subgroup of
finite index isomorphic to $\IZ^2$, which is only possible in the
case~\ref{the:virtually_poly_Z-groups:ge_2} and not in the
case~\ref{the:virtually_poly_Z-groups:1_and_N_GV_is_G} under consideration.

Next, we want to show that $\vcd(G) \ge 4$ in the
case~\ref{the:virtually_poly_Z-groups:1_and_N_GV_is_G:vcd-1} and $\vcd(G)\ge 3$
in the case~\ref{the:virtually_poly_Z-groups:1_and_N_GV_is_G:vcd}. We can
assume that $C_0\subseteq G$ is central and that $G$ is poly-$\IZ$ for,
otherwise, we can replace $G$ by a subgroup of finite index.  Choose an element
$g \in G$ whose image under the projection $G \to G/C_0$ is an element of
infinite order.  Let $D$ be the cyclic subgroup generated by $g \in G$. Then $D
\cap C_0 = \{1\}$ and $N_GD$ contains a subgroup isomorphic to $\IZ^2$. Hence
$\vcd(N_GD) \ge 2$ and $[D] \not= [C_0]$. The assumptions imply $\vcd(D) \le
\vcd(G)-2$, which in turn implies $\vcd(G) \ge 4$ in the
case~\ref{the:virtually_poly_Z-groups:1_and_N_GV_is_G:vcd-1} and $\vcd(G) \ge
3$ in the case~\ref{the:virtually_poly_Z-groups:1_and_N_GV_is_G:vcd}.

In order to get a model for $\edub{G}$ of the desired dimension, we have to
investigate the $G$-pushout~\eqref{basic_G-pushout} more closely and use some
machinery. Let $\edub{G}$ be the $G$-$CW$-model obtained from the
$G$-pushout~\eqref{basic_G-pushout}. Given a natural number $d \ge 3$, we will
show that there exists a $G$-$CW$-complex of dimension $d$ which is
$G$-homotopy equivalent to $\edub{G}$ by checking that for every contravariant
$\IZ\Or(G)$-module $M$ the Bredon cohomology $H^{d+1}_{\IZ\Or(G)}(\edub{G};M)$
vanishes and for every $i \ge d+1$ the homology of $\edub{G}^H$ vanishes for $H
\subseteq G$ (see~\cite[Proposition~11.10 on page~221 and Theorem~13.19 on
page~268]{Lueck(1989)}).  The latter condition is obviously satisfied since
$\edub{G}^H$ is contractible or empty.  It remains to prove
$H^{d+1}_{\IZ\Or(G)}(\edub{G};M) = 0$.

Let $C \in I$ be the unique element for which $\vcd(N_GC) = \vcd(G)$. After
replacing $C$ by $C_0$, we can actually assume that $N_GC=G$. Let us put $d =
\vcd(G)$ if $\vcd(N_GD) = \vcd(G) -1$ for at least one $D \in I$ with $D \not=
C$ and $d = \vcd(G) - 1$ if $\vcd(N_GD) \le \vcd(G) -2$ holds for all $D \in I$
with $D \not= C$. Then we have $d\ge 3$ and, in all cases, $\vcd(N_GD) < d$ for
all $D \in I$ with $D \not= C$. In particular, $H^{k}_{\IZ\Or(G)}(G
\times_{N_GD} \eub{N_GD};M)$ as well as $H^{k}_{\IZ\Or(G)}(G \times_{N_GD}
\eub{W_GD};M)$ vanish for all $D \in I$ with $D \not= C$ if $k \ge d$. The long
exact Mayer-Vietoris sequence for the Bredon cohomology associated to the
$G$-pushout~\eqref{basic_G-pushout} yields the exact sequence
\begin{multline*}
  H^{d}_{\IZ\Or(G)}(\eub{G};M) \xrightarrow{\id} H^{d}_{\IZ\Or(G)}(\eub{G};M) \to
  H^{d+1}_{\IZ\Or(G)}(\edub{G};M) \\ \to H^{d+1}_{\IZ\Or(G)}(\eub{G};M) \xrightarrow{\id}
  H^{d+1}_{\IZ\Or(G)}(\eub{G};M).
\end{multline*}
Hence, $H^{d+1}_{\IZ\Or(G)}(\edub{G};M)$ is trivial and, as explained above, we
can find a $d$-dimensional $G$-$CW$-model for $\edub{G}$.  We conclude that
$\hdim^G(\edub{G}) \le \vcd(G) -1$ in the
case~\ref{the:virtually_poly_Z-groups:1_and_N_GV_is_G:vcd-1} and
$\hdim^G(\edub{G}) \le \vcd(G)$ in the
case~\ref{the:virtually_poly_Z-groups:1_and_N_GV_is_G:vcd}.

Finally, we prove $\hdim^G(\edub{G}) \ge \vcd(G) -1$ in the
case~\ref{the:virtually_poly_Z-groups:1_and_N_GV_is_G:vcd-1} as well as
$\hdim^G(\edub{G}) \ge \vcd(G)$ in the
case~\ref{the:virtually_poly_Z-groups:1_and_N_GV_is_G:vcd}.  The first
inequality follows from Corollary~\ref{cor:dim(eub(G)_le_dim(edub(G)_plus_1}.
Concerning the second, we can assume that $G$ is poly-$\IZ$. Consider the long
exact Mayer-Vietoris sequence associated to the pushout
\eqref{quotient_of_basic_G-pushout}. Since a poly-$\IZ$-group is torsionfree,
$N_GC = G$ and $\cd(EW_GD) \le d-2$ for $D \in I$ with $D \not= C$, 
it yields the long exact sequence
\begin{multline*}
  H^{d-1}(G) \xrightarrow{\id \times j} H^{d-1}(G) \times \prod_{D \in I, D\not= C}
  H^{d-1}(N_GD) \to H^d(G\backslash\edub{G}) \\
  \to H^{d}(G) \xrightarrow{\id} H^{d}(G).
\end{multline*}
As $H^{d-1}(N_GD)$ is non-trivial if $D \in I$ is such that $\vcd(N_GD)= d-1$
(see
Lemma~\ref{lem:prop_virtually_poly_Z-gr}~\ref{lem:prop_virtually_poly_Z-gr:top_homology}),
we conclude $H^d(G\backslash \edub{G}) \not= 0$ and hence $\hdim^G(\edub{G})\ge
d=\vcd(G)$. This finishes the proof of
Theorem~\ref{the:virtually_poly_Z-groups}.
\end{proof}

\begin{corollary} \label{cor:hdim(edub(G)_finite_quit_or_index}
Let $G$ be a virtually poly-$\IZ$ group. Suppose that the group $G'$ is a
subgroup of $G$ of finite index or that $G'$ is the quotient of $G$ by a finite
normal subgroup. Then
\begin{align*}
    \hdim^{G'}(\eub{G'}) &= \hdim^G(\eub{G}); \\
    \hdim^{G'}(\edub{G'}) &= \hdim^G(\edub{G}).
\end{align*}
\end{corollary}
\begin{proof}
We get $\vcd(G') = \vcd(G)$ from
Lemma~\ref{lem:prop_virtually_poly_Z-gr}~\ref{lem:prop_virtually_poly_Z-gr:exact_sequences_and_vcd}.
Hence, $\hdim^{G'}(\eub{G'}) = \hdim^G(\eub{G})$ follows immediately from
Theorem~\ref{the:virtually_poly_Z-groups}.

The second assertion of the corollary will follow from
assertions~\ref{the:virtually_poly_Z-groups:empty},~\ref{the:virtually_poly_Z-groups:1_and_N_GV_is_G}
and~\ref{the:virtually_poly_Z-groups:ge_2} of
Theorem~\ref{the:virtually_poly_Z-groups}. We only need to show that every
subgroup of $G'$ of finite index has a finite center if and only if every
subgroup of $G$ of finite index has a finite center and that $G'$ contains a
finite index subgroup whose center contains $\IZ^2$ if and only if $G$ contains
a finite index subgroup whose center contains $\IZ^2$.

We begin with the case that $G'$ is a subgroup of finite index in $G$.  A
finitely generated abelian group $A$ satisfies $\vcd(A) \ge 2$ if and only if
it contains a subgroup isomorphic to $\IZ^2$. Obviously, a subgroup $H
\subseteq G$ satisfies $[G:H] < \infty$ if and only if $[G':(G'\cap H)] <
\infty$. Now the claims follow from the observation that for a subgroup
$H\subseteq G$ we have $\cent(H) \cap G' \subseteq \cent(H \cap G')$ and
$[\cent(H): (\cent(H) \cap G')] < \infty$, which implies that $\vcd(\cent(H))
\le \vcd(\cent(H \cap G'))$ holds.

It remains to treat the case that there is an extension $1 \to F \to G
\xrightarrow{p} G' \to 1$ in which the group $F$ is finite. Let $K \subseteq G$
be a subgroup, then $p(\cent(K)) \subseteq \cent(p(K))$, which implies
\begin{equation} \label{cor:hdim(edub(G)_finite_quit_or_index:eq1}
    \vcd(\cent(K)) \le \vcd(\cent(p(K))).
\end{equation}
Let $H \subseteq G'$ be a subgroup. Since $G$ is virtually poly-$\IZ$, by
Lemma~\ref{lem:prop_virtually_poly_Z-gr}~\ref{lem:prop_virtually_poly_Z-gr:inheritance}
the same is true for the subgroup $p^{-1}(H)$. In particular, $p^{-1}(H)$ is
finitely generated, say by the finite set $\{s_1, s_2, \ldots, s_r\}$. Consider
the map
\begin{equation*}
    \gamma \colon p^{-1}(\cent(H)) \to \prod_{i=1}^r F, \quad h\mapsto \{[h,s_i]
    \mid i = 1,2, \ldots r\}.
\end{equation*}
We have $p(\cent(p^{-1}(H)) \subseteq \cent(p(p^{-1}(H)) = \cent(H)$ and hence
$\cent(p^{-1}(H))$ is a subgroup of $p^{-1}(\cent(H))$. Now one easily checks
that $\gamma$ induces an injection $p^{-1}(\cent(H))/\cent(p^{-1}(H)) \to
\prod_{i=1}^r F$, which means that $\cent(p^{-1}(H))$ has finite index in
$p^{-1}(\cent(H))$. As $F$ is finite, we conclude that
\begin{equation} \label{cor:hdim(edub(G)_finite_quit_or_index:eq2}
    \vcd(\cent(p^{-1}(H))) = \vcd(\cent(H)).
\end{equation}
Obviously, $[G:K] < \infty$ if and only if $[G':p(K)] < \infty$, and $[G':H] <
\infty$ if and only if $[G:p^{-1}(H)] < \infty$. Thus, the claim follows
from~\eqref{cor:hdim(edub(G)_finite_quit_or_index:eq1}
and~\eqref{cor:hdim(edub(G)_finite_quit_or_index:eq2}.
\end{proof}


\subsection{Extensions of free abelian groups}
\label{subsec:Extensions_of_free_abelian_groups}

The following is a consequence of
Theorem~\ref{the:virtually_poly_Z-groups}~\ref{the:virtually_poly_Z-groups:ge_2}
and Corollary~\ref{cor:hdim(edub(G)_finite_quit_or_index}. It has already been
proved for crystallographic groups by
Connolly-Fehrmann-Hartglass~\cite{Connolly-Fehrmann-Hartglass(2006)}.

\begin{example}[Virtually $\IZ^n$-groups]
\label{exa:mindim(G)_G_virtually_Zn}
Let $G$ be a group which contains $\IZ^n$ as subgroup of finite index. If $n
\ge 2$, then
\begin{equation*}
    \hdim^G(\edub{G}) = n+1.
\end{equation*}
If $ n \le 1$, then $\hdim^G(\edub{G}) = 0$.
\end{example}

\begin{theorem} \label{the:Extensions_1_to_Zn_to_G_to_Z_to_1}
Let $1 \to \IZ^n \to G \xrightarrow{p} \IZ \to 1$ be a group extension for $n
\ge 2$.  Let $\phi \colon \IZ^n \to \IZ^n$ be the automorphism given by
conjugation with some $g \in G$ which is mapped to a generator of $\IZ$ under
$p$.

Then there is the following dichotomy:

\begin{enumerate}

\item \label{the:Extensions_1_to_Zn_to_G_to_Z_to_1_le_1}
We have $\cd(\ker(\phi^k - \id)) \le 1$ for all $k \in \IZ$.  In this case,
\begin{equation*}
    \hdim^G(\edub{G}) = n+1;
\end{equation*}

\item \label{the:Extensions_1_to_Zn_to_G_to_Z_to_1_ge_2}
We have $\cd(\ker(\phi^k - \id)) \ge 2$ for some $k \in \IZ$.  In this case,
\begin{equation*}
    \hdim^G(\edub{G}) = n+2.
\end{equation*}

\end{enumerate}
\end{theorem}
\begin{proof}
We first deal with the case that $\phi$ is periodic, i.e., there exists $k \ge
1$ such that $\phi^k = \id$. Since $n \ge 2$, we are then in the situation of
case~\ref{the:Extensions_1_to_Zn_to_G_to_Z_to_1_ge_2}. The group $G$ contains a
subgroup of finite index which is isomorphic to $\IZ^{n+1}$ and, hence,
$\hdim^G(\edub{G}) = n+2$ by Example~\ref{exa:mindim(G)_G_virtually_Zn}. Thus,
we can assume in the remainder that $\phi$ is not periodic.

Next, we show that the epimorphism $p \colon G \to \IZ$ satisfies $p(V)=0$ for
all $[V]\cdot G\in[\calvcyc\setminus\calfin]_f/G$. Suppose to the contrary that
there exists $[V] \in [\calvcyc\setminus \calfin]$ such that $N_G[V]$ has
finite index in $G$ and $p(V) \not= 0$. Because of
Lemma~\ref{lem:N_G[V]_is_N_GV}, we can assume without loss of generality that
$V$ is infinite cyclic and $N_G[V] = N_GV$. Since $C_GV$ has finite index in
$N_GV$, we conclude that $C_GV$ has finite index in $G$.  Let $v \in V$ be a
generator of $V$. Then $\phi^{p(v)} \colon \IZ^n \to \IZ^n$ is given by
conjugation with $v$ and, hence, is the identity on the subgroup $\IZ^n \cap
C_GV$ of finite index in $\IZ^n$. This shows that $\phi^{p(v)} = \id$,
contradicting the assumption that $\phi$ is not periodic.

Now, consider $[V] \cdot G\in [\calvcyc\setminus\calfin]_f/G$. Then $p(V) = 0$
implies $\IZ^n \subseteq N_G[V]$ and $N_G[V] = p^{-1}(p(N_G[V]))$. One easily
checks that
\begin{equation*}
    p(N_G[V]) = \{k \in \IZ \mid V \subseteq \ker(\phi^k -\id) \;\text{or } V
    \subseteq \ker(\phi^k +\id) \}.
\end{equation*}%
\ref{the:Extensions_1_to_Zn_to_G_to_Z_to_1_le_1} Obviously, $\ker(\phi^k - \id)
= 0$ holds for all $k \in \IZ$ if and only if $[\calvcyc\setminus\calfin]_f/G$
is empty. In this case, $\hdim^G(\edub{G}) = n+1$ because of
Theorem~\ref{the:virtually_poly_Z-groups}~\ref{the:virtually_poly_Z-groups:empty}.
We are left to treat the case that there exists $n \ge 1$ such that
$\cd(\ker(\phi^n - \id)) = 1$, while $\cd(\ker(\phi^k - \id)) \le 1$ holds for
all $k \ge 1$. Put $V := \ker(\phi^n - \id)$. Then the the index of $N_GV$ in
$G$ is finite, and so $[V] \cdot G$ belongs to
$[\calvcyc\setminus\calfin]_f/G$. Consider another element $[W] \cdot
G\in[\calvcyc\setminus\calfin]_f/G$. There exists $k \ge 1$ such that $W
\subseteq \ker(\phi^k -\id)$ or $W \subseteq \ker(\phi^k +\id)$, hence $W
\subseteq \ker(\phi^{2kn} -\id)$. Since also $V \subseteq \ker(\phi^{2kn}
-\id)$ and $\cd(\ker(\phi^{2kn} -\id)) \le 1$, we conclude $[V] \cdot G = [W]
\cdot G$, and we have shown that $[\calvcyc\setminus\calfin]_f/G=\{[V]\cdot
G\}$ contains precisely one element $[V] \cdot G$.

Now $\hdim^G(\edub{G}) = n+1$ follows from
Theorem~\ref{the:virtually_poly_Z-groups}~\ref{the:virtually_poly_Z-groups:1_and_N_GV_is_G:vcd}.
Namely, any infinite cyclic subgroup $C\subseteq\IZ^n$ such that $C\cap
V=\{1\}$ yields an element $[C]\cdot G\in[\calvcyc\setminus\calfin]_f/G$ such
that $\cd(N_G[C])=n=\cd(G)-1$.
\\[1mm]\ref{the:Extensions_1_to_Zn_to_G_to_Z_to_1_ge_2}
Choose $k \ge 1$ with $\cd(\ker(\phi^k - \id)) \ge 2$ and put
$G':=p^{-1}(k\IZ)$. Then $G'$ is a subgroup of $G$ of finite index whose center
contains $\IZ^2$, and $\hdim^G(\edub{G}) = n+2$ follows from
Theorem~\ref{the:virtually_poly_Z-groups}~\ref{the:virtually_poly_Z-groups:ge_2}.
\end{proof}

\begin{theorem} \label{the:Central_extensions_of_virtually_Zn-groups}
Let $ 1 \to \IZ^m \to G \xrightarrow{p} \IZ^n \to 1$ be a group extension such
that $\cent(G) = \IZ^m$. Then
\begin{equation*}
    \hdim^G(\eub{G}) = \cd(G) = m + n
\end{equation*}
and
\begin{equation*}
    \hdim^G(\edub{G}) =
    \begin{cases}
        0 & \text{if } m + n \le 1; \\
        n+1 & \text{if } m = 1 \; \text{and } n \ge 1; \\
        m+n+1 & \text{if } m \ge 2 \; \text{and } n \ge 1. \\
    \end{cases}
\end{equation*}
\end{theorem}
\begin{proof}
The group $G$ is poly-$\IZ$ by
Lemma~\ref{lem:prop_virtually_poly_Z-gr}~\ref{lem:prop_virtually_poly_Z-gr:inheritance}.
We get $\hdim^G(\eub{G})=\cd(G)=m+n$ from
Theorem~\ref{the:virtually_poly_Z-groups} and
Lemma~\ref{lem:prop_virtually_poly_Z-gr}~\ref{lem:prop_virtually_poly_Z-gr:exact_sequences_and_vcd}.

If $m + n \le 1$, then $G$ is infinite cyclic, so $\hdim^G(\edub{G}) = 0$. If $
m = 0$, then $G = \IZ^n$ and the claim follows from
Example~\ref{exa:mindim(G)_G_virtually_Zn}. If $ m\ge 2$, the claim follows
from
Theorem~\ref{the:virtually_poly_Z-groups}~\ref{the:virtually_poly_Z-groups:ge_2}.
It remains to treat the case $m = 1$ and $n \ge 1$, in which we want to show
for an infinite cyclic subgroup $C \subseteq G$ with $[C]\not=[\cent(G)]$ that
$\cd(N_G[C])=n$.

Consider the canonical projection $p\colon G \to G/\cent(G) \cong \IZ^n$.
According to Lemma~\ref{lem:N_G[V]_is_N_GV}, we can assume that $N_GC =
N_G[C]$. Since $C_GC$ has finite index in $N_GC$, it suffices to prove that
$\cd(C_GC) = n$. Fix a generator $z\in C$ and an element $g\in G$. Let $c(g)
\colon G \to G$ be the automorphism sending $g'$ to $gg'g^{-1}$, then $p\circ
c(g)(z)=p(z)$. There is, therefore, precisely one element $\sigma(g) \in
\cent(G)$ satisfying $c(g)(z)=z\sigma(g)$. A straightforward computation for
$g_1,g_2 \in G$ shows that $z\sigma(g_1g_2)=z\sigma(g_1)\sigma(g_2)$, which
implies $\sigma(g_1g_2)=\sigma(g_1)\sigma(g_2)$. Moreover, $\sigma(1) = 1$.
Thus, we have defined a group homomorphism $\sigma\colon G\to\cent(G)$ whose
kernel is $C_GC$. As $[C]\not=[\cent(G)]$ by assumption, $\cent(G)$ cannot
contain $C$, so $\ker(\sigma)\not=G$. This means that $\sigma$ is non-trivial
and that, hence, its image is an infinite cyclic group. We conclude from
Lemma~\ref{lem:prop_virtually_poly_Z-gr}~\ref{lem:prop_virtually_poly_Z-gr:exact_sequences_and_vcd}
that $\cd(C_GC) = n$, as we wanted to show.

The above implies that we are in the situation of
Theorem~\ref{the:virtually_poly_Z-groups}~\ref{the:virtually_poly_Z-groups:1_and_N_GV_is_G:vcd},
and we get $\hdim^G(\edub{G}) = \cd(G) = n+1$.
\end{proof}

\begin{theorem} \label{the:Hei_rtimesZ}
Let $1 \to \IZ \to H \xrightarrow{p} \IZ^n \to 1$ be a group extension such
that $\IZ = \cent(H)$ and $n \ge 2$.  Consider a group automorphism $f \colon H
\to H$. It sends $\cent(H)$ to $\cent(H)$. Let $\overline{f} \colon \IZ^n \to
\IZ^n$ be the group automorphism induced by $f$, i.e., which satisfies $p \circ
f = \overline{f} \circ p$. Then $G:=H\rtimes_f\IZ$ satisfies $\cd(G)=n+2$, we
have
\begin{equation*}
    \hdim^G(\eub{G}) = n+2,
\end{equation*}
and precisely one of the following cases occurs:

\begin{enumerate}

\item \label{the:Hei_rtimesZ:n_plus_1}
We have $ker(\overline{f}^k - \id) = 0$ for every $k \in \IZ$, $k \not= 0$. In
this case,
\begin{equation*}
    \hdim^G(\edub{G}) = n+1;
\end{equation*}

\item \label{the:Hei_rtimesZ:n_plus_2}
We have $\ker(\overline{f}^k - \id) \not= 0$ for some $k \in \IZ$, $k \not= 0$
and $\overline{f}$ is not periodic, i.e., there is no $l\in\IZ$, $l\ge 1$ such
that $\overline{f}^l = \id$.  In this case,
\begin{equation*}
    \hdim^G(\edub{G}) = n+2;
\end{equation*}

\item \label{the:Hei_rtimesZ:n_plus_3}
The map $\overline{f}$ is periodic. In this case,
\begin{equation*}
    \hdim^G(\edub{G}) = n+3.
\end{equation*}

\end{enumerate}
\end{theorem}
\begin{proof}
Let $\pr \colon G = H \rtimes_f \IZ \to \IZ$ be the obvious projection, then
$G$ fits into an exact sequence $1\to H\to G\xrightarrow{\pr}\IZ\to 1$.
Therefore, $\hdim^G(\eub{G})=\cd(G)=n+2$ follows from
Lemma~\ref{lem:prop_virtually_poly_Z-gr}~\ref{lem:prop_virtually_poly_Z-gr:inheritance}
and~\ref{lem:prop_virtually_poly_Z-gr:exact_sequences_and_vcd} and
Theorem~\ref{the:virtually_poly_Z-groups}.

Let $V \subseteq G$ be an infinite cyclic subgroup.  Then $V \subseteq
\cent(H)$ if and only if $[V] = [\cent(H)]$ and, in this case, $N_G[V] = G$.
Provided that $[V] \not= [\cent(H)]$, we first want to show
\begin{equation} \label{the:Hei_rtimesZ:(1)}
    \cd(N_G[V])=
    \begin{cases}
        n & \text{if } \pr(N_G[V]) = 0; \\
        n+1 & \text{if } \pr(N_G[V]) \not=  0 \\
        & \phantom{if } \text{and } \pr(V) = 0; \\
        2+\cd(p(N_G[V] \cap H)) & \text{if } \pr(V)\not= 0.
    \end{cases}
\end{equation}

Suppose that $\pr(V) = 0$, i.e., $V \subseteq H$.  We have already seen in the
proof of Theorem~\ref{the:Central_extensions_of_virtually_Zn-groups} that
$\cd(N_{H}[V]) = n$ for $[V] \not= [\cent(H)]$. It follows from the short exact
sequence
\begin{equation} \label{the:Hei_rtimesZ:(4)}
    1 \to N_{H}[V] = N_G[V] \cap H \to N_G[V] \to \pr(N_G[V]) \to 1
\end{equation}
and
Lemma~\ref{lem:prop_virtually_poly_Z-gr}~\ref{lem:prop_virtually_poly_Z-gr:exact_sequences_and_vcd}
that $\cd(N_G[V]) = n + \cd(\pr(N_G[V]))$. If $\pr(N_G[V]) = 0$, then
$\cd(\pr(N_G[V])) =0$. If $\pr(N_G[V]) \not= 0$, then $\pr(N_G[V]) \cong \IZ$,
which implies $\cd(\pr(N_G[V])) = 1$.

Now suppose that $\pr(V) \not= 0$.  Then $\pr(V)$ is infinite cyclic, so the
same must be true for $\pr(N_G[V])$. Let $w$ be a generator of $\cent(H)$ and
$v$ be a generator of $V$. Since $f$ induces $\pm \id$ on $\cent(H)$, we have
$vwv^{-1} = w^{\pm 1}$. This implies $vw^2v^{-1} = w^2$, thus $w^2 \in C_GV
\subseteq N_G[V]$. Therefore, $N_G[V]\cap\cent(H)$ is an infinite cyclic group,
and the short exact sequences
\begin{equation*}
    1 \to N_G[V] \cap \cent(H) \to N_G[V] \cap H \to p(N_G[V] \cap H) \to 1
\end{equation*}
and \eqref{the:Hei_rtimesZ:(4)} together with
Lemma~\ref{lem:prop_virtually_poly_Z-gr}~\ref{lem:prop_virtually_poly_Z-gr:exact_sequences_and_vcd}
imply $\cd(N_G[V]) = 2 + \cd(p(N_G[V] \cap H))$.  This finishes the proof
of~\eqref{the:Hei_rtimesZ:(1)}.

Next, we show in the case $\pr(V) \not= 0$ that
\begin{equation} \label{the:Hei_rtimesZ:(2)}
    p(N_G[V] \cap H) \subseteq \ker(\overline{f}^k - \id)
    \quad \text{for some } k \in \IZ, k \not= 0.
\end{equation}
Because of Lemma~\ref{lem:N_G[V]_is_N_GV}, we can assume without loss of
generality that $N_G[V] = N_GV$.  Recall that $C_GV$ has finite index in
$N_GV$. Hence it suffices to show that $p(C_GV\cap H) \subseteq
\ker(\overline{f}^k - \id)$ if $k \in \pr(V)$. Let $t \in G = H \rtimes_f \IZ$
be given by the generator of $\IZ$. Then conjugation with $t$ induces the
automorphism $f$ on $H$, and $ht^k \in V$ for some $h \in H$. We get for $u \in
C_GV\cap H$
\begin{multline*}
    \overline{f}^k(p(u)) = p \circ f^k(u) = p(t^kut^{-k})
    = p(h) p(t^kut^{-k}) p(h)^{-1} \\
    = p(ht^ku(ht^k)^{-1}) = p(u),
\end{multline*}
which finishes the proof of~\eqref{the:Hei_rtimesZ:(2)}.

Finally, we show in the case $\pr(N_G[V]) \not= 0$ and $\pr(V) = 0$ that
\begin{equation} \label{the:Hei_rtimesZ:(3)}
    p(V) \subseteq \ker(\overline{f}^k - \id)
    \quad \text{for some } k \in \IZ, k \not= 0.
\end{equation}
Again, we can assume that $N_G[V] = N_GV$.  Since $C_GV$ has finite index in
$N_GV$, we have $\pr(C_G[V]) \not= 0$, so we can choose a non-trivial element
$k \in \pr(C_G[V])$. Then $ht^k \in C_GV$ for some $h \in H$, and the same
calculation as above shows that $\overline{f}^k(p(v)) = p(v)$ holds for $v\in
V$. This finishes the proof of~\eqref{the:Hei_rtimesZ:(3)}.  Now we are ready
to prove the assertions appearing in Theorem~\ref{the:Hei_rtimesZ}.
\\[1mm]\ref{the:Hei_rtimesZ:n_plus_1}
Consider an infinite cyclic subgroup $V \subseteq G$ such that $[V] \not=
[\cent(H)]$. We first deal with the case $\pr(V) = 0$. Then 
$\pr(N_G[V]) = 0$ since, otherwise,
\eqref{the:Hei_rtimesZ:(3)} would yield $V \subseteq \cent(H)$. Using
\eqref{the:Hei_rtimesZ:(1)}, the conclusion is that $\cd(N_G[V])=n$. 
Since $\cd(G) = n+2$, we get
$\hdim^G(\edub{G}) = n+1$ from
Theorem~\ref{the:virtually_poly_Z-groups}~\ref{the:virtually_poly_Z-groups:1_and_N_GV_is_G:vcd-1}.
It remains to consider the case $\pr(V) \not= 0$. Then we conclude
$\cd(N_G[V]) = 2 + \cd(p(N_GV \cap H))$ from~\eqref{the:Hei_rtimesZ:(1)}.
Because of~\eqref{the:Hei_rtimesZ:(2)} we get $\cd(N_G[V]) = 2$.
Since $n \ge 2$, we have $\vcd(G) \ge 4$. Now $\hdim^G(\edub{G}) = n+1$ follows from
Theorem~\ref{the:virtually_poly_Z-groups}%
~\ref{the:virtually_poly_Z-groups:1_and_N_GV_is_G:vcd-1}.
\\[2mm]\ref{the:Hei_rtimesZ:n_plus_2}
We claim that there exists an infinite cyclic subgroup $V \subseteq H$ such
that $[V] \not= [\cent(H)]$, $\pr(V)=0$ and $\pr(N_GV) \not= 0$. To do so, we
choose an even $k \in \IZ$ such that $\ker(\overline{f}^k - \id) \not= 0$, an
element $h_0 \in H$ such that $0\not=p(h_0) \in \ker(\overline{f}^k - \id)$,
and define $V$ to be the infinite cyclic group generated by $h_0$. Then
$p(f^k(h_0)) = \overline{f}^k \circ p(h_0) = p(h_0)$. Hence, we can find $z_0
\in \cent(H) = \IZ$ such that $t^kh_0t^{-k} = f^k(h_0) = h_0z_0$. Since $k$ is
even, $t^kz_0t^{-k} = z_0$. As $h_0 \notin \cent(H)$, we can find $h_1 \in H$
such that $h_1h_0h_1^{-1} \not= h_0$. However, $p(h_1h_0h_1^{-1}) =
p(h_1)p(h_0)p(h_1^{-1}) = p(h_0)$ and so $h_1h_0h_1^{-1} = h_0z_1$ for some
non-trivial $z_1 \in \cent(H)$. Choose integers $m_0\not=0$ and $m_1$ such that
$z_0^{m_0}= z_1^{m_1}$, and put $u = h_1^{-m_1}t^{m_0k}$. From
$t^{m_0k}h_0t^{-m_0k} = h_0z_0^{m_0}$ and
$h_1^{-m_1}h_0h_1^{m_1}=h_0z_1^{-m_1}$ we then conclude $uh_0u^{-1}=h_0$ and
$\pr(u) = km_0$. This finishes the proof of the above claim. It follows from
\eqref{the:Hei_rtimesZ:(1)} that $\cd(N_G[V]) = n+1$.

Since $\cd(\ker(\overline{f}^k - \id))= n$ implies $\overline{f}^k = \id$ but
$\overline{f}$ is not periodic, we have $\cd(\ker(\overline{f}^k - \id)) \le
n-1$ for every $k \in \IZ$ with $k \not= 0$. Thus, \eqref{the:Hei_rtimesZ:(1)}
and~\eqref{the:Hei_rtimesZ:(2)} show that $\cd(N_G[V]) \le n+1$ holds for all
infinite cyclic subgroups $V\subseteq G$ satisfying $[V] \not= [\cent(H)]$.
Theorem~\ref{the:virtually_poly_Z-groups}~\ref{the:virtually_poly_Z-groups:1_and_N_GV_is_G:vcd}
now implies $\hdim^G(\edub{G}) = n+2$.
\\[2mm]\ref{the:Hei_rtimesZ:n_plus_3}
We want to show the existence of a finite index subgroup of $G$ whose center
contains $\IZ^2$. Then $\hdim^G(\edub{G}) = n+3$ follows from
Theorem~\ref{the:virtually_poly_Z-groups}~\ref{the:virtually_poly_Z-groups:ge_2}.
Let us begin by choosing an even $k\in\IZ$ such that $\overline{f}^k = \id$. We
define a group homomorphism $\sigma_{f^k} \colon H \to \cent(H) = \IZ$ by
$\sigma_{f^k}(h):=h^{-1}f^k(h)$ for $h\in H$. Since $k$ is even, the
restriction of $f^k$ to $\cent(H)=\IZ$ is the identity, which implies that the
restriction of $\sigma_{f^k}$ to $\cent(H)$ is trivial. Hence $\sigma_{f^k}$
factorizes through the projection $p\colon H \to H/\cent(H) = \IZ^n$, yielding
homomorphism $\overline{\sigma}_{f^k} \colon \IZ^n \to \IZ$. We have $f^k(h) =
h \overline{\sigma}_{f^k}(p(h))$ for $h \in H$.

Conjugation with $h_0 \in H$ induces a group automorphism $c(h_0) \colon H \to
H$ whose restriction to $\cent(H)$ is the identity and for which the induced
automorphism $\overline{c(h_0)}$ of $H/\cent(H) = \IZ^n$ is the identity. Just
as above for $f^k$, we define a homomorphism $\overline{\sigma}_{c(h_0)} \colon
\IZ^n \to \IZ$ such that $c(h_0)(h) = h_0hh_0^{-1} = h
\overline{\sigma}_{c(h_0)}(p(h))$ holds for $h \in H$, and one easily checks
that
\begin{equation*}
    \overline{\sigma}_{c(h_0h_1)} = \overline{\sigma}_{c(h_0)\circ c(h_1)}
    = \overline{\sigma}_{c(h_0)} + \overline{\sigma}_{c(h_0)}
\end{equation*}
for $h_0, h_1 \in H$. That way, we obtain an injective group homomorphism
\begin{equation*}
    \overline{\sigma} \colon H/\cent(H) \to \hom_{\IZ}(H/\cent(H),\IZ),
    \quad \overline{h} \mapsto \overline{\sigma}_{c(h)}.
\end{equation*}
Since its source and target are isomorphic to $\IZ^n$, its image has finite
index. As $\overline{\sigma}_{f^{lk}} = l \cdot \overline{\sigma}_{f^k}$ holds
for all $l \in \IZ$, it follows that we can find an integer $l \ge 1$ and $h
\in H$ such that $\overline{\sigma}_{f^{lk}} = \overline{\sigma}_{c(h)}$, which
implies $f^{kl} = c(h)$.

The subgroup $\pr^{-1}(kl \cdot \IZ) \subseteq G$ has finite index in $G$ and
is isomorphic to the semidirect product $H \rtimes_{f^{kl}} \IZ$ with respect
to the automorphism $f^{kl} \colon H \to H$. Since $f^{kl}$ is an inner
automorphism, $\pr^{-1}(kl \cdot \IZ)$ is isomorphic to the direct product $H
\times \IZ$ and, hence, contains $\cent(H) \times \IZ \cong \IZ^2$ in its
center.  This finishes the proof of Theorem~\ref{the:Hei_rtimesZ}.
\end{proof}

\begin{example}[Three-dimensional Heisenberg group]
\label{the:Heisenberg_group_extension_by_Z}
The three-di\-men\-sional Heisenberg group is the group given by the
presentation $\Hei = \langle u,v,z \mid [u,v] = z, [u,z] = 1, [v,z] = 1
\rangle$. Define automorphisms $f_k \colon \Hei \to \Hei$ for $k = -1,0,1$ by
\begin{align*}
    f_{-1}(u) &= u^3v, & f_{-1}(v) &= u^2v, & f_{-1}(z) &= z; \\
    f_0(u) &= u, & f_0(v) &= uv, & f_0(z) &= z; \\
    f_1(u) &= u, & f_1(v) &= v, & f_1(z) &= z.
\end{align*}
Put $G= \Hei \rtimes_{f_k} \IZ$ for $k = -1,0,1$. Then
Theorem~\ref{the:Hei_rtimesZ} implies
\begin{align*}
    \hdim^{G_k}(\eub{G_k)} &= 4; \\
    \hdim^{G_k}(\edub{G_k)} &= 4 + k.
\end{align*}
\end{example}


\subsection{Models for $\EGF{G}{\calsfg}$}
\label{subsec:Models_for_EGcalsfg}

Let $\calsfg$ be the family of subgroups of $G$ which are contained in some
finitely generated subgroup, i.e., $H\in\calsfg$ if and only if there exists a
finitely generated subgroup $H'\subseteq G$ such that $H\subseteq H'$. This is
the smallest family of subgroups of $G$ which contains all finitely generated
subgroups of $G$.

\begin{example} \label{exa:calf_is_fin_gen_sub}
Suppose that $G$ is countable. Then there exists a one-di\-men\-sional model
for $\EGF{G}{\calsfg}$ because of Theorem~\ref{the:model_for_colimits} and
Lemma~\ref{lem:hdim(I)}.  In particular, there exists a one-dimensional model
for $\eub{G}$ or $\edub{G}$ if $G$ is locally finite or locally virtually
cyclic, respectively.
\end{example}

\begin{theorem} \label{thm:locally_F_model}
Let $G$ be an infinite group of cardinality $\aleph_n$. Let $\calsfg$ be the
family of subgroups defined above. Then
\begin{equation*}
    \hdim^G(\EGF{G}{\calsfg}) \le n+1.
\end{equation*}
\end{theorem}
\begin{proof}
We use induction on $n\in\IN$. The case $n=0$ has already been settled in
Example~\ref{exa:calf_is_fin_gen_sub}. The induction step from $n-1$ to $n\ge
1$ is done as follows.

We can write $G=\bigcup_{\alpha<\omega_n}G_\alpha$ for subgroups $G_\alpha$ of
cardinality $\aleph_{n-1}$ such that $G_\alpha\subseteq G_\beta$ if
$\alpha\le\beta$. By induction hypothesis, for $\alpha<\omega_n$ there are
models $X_{\alpha}$ for $\EGF{G_{\alpha}}{\calsfg}$ of dimension $n$. The
induction step now must provide us with an $(n+1)$-dimensional model for
$\EGF{G}{\calsfg}$. If we set $G_{\omega_n} := G$, this will be accomplished by
using transfinite induction for $\alpha\le\omega_n$ to construct
$(n+1)$-dimensional models $Y_{\alpha}$ for $\EGF{G_{\alpha}}{\calsfg}$ such
that $G_{\gamma}\times_{G_{\alpha}}Y_{\alpha} \subseteq G_{\gamma}
\times_{G_{\beta}} Y_{\beta}$ is a $G_{\gamma}$-subcomplex if
$\alpha\le\beta\le\gamma$.

Let $Y_0 := X_0$. Now suppose that $\alpha$ has got a predecessor. Then the
universal property of $\EGF{G_{\alpha}}{\calsfg}$ yields a $G_{\alpha}$-map
$f_{\alpha}\colon G_{\alpha}\times_{G_{\alpha-1}} X_{\alpha-1} \to X_{\alpha}$,
which we can assume to be cellular by the equivariant cellular approximation
theorem. We define $Y_{\alpha}$ by the $G_{\alpha}$-pushout
\begin{equation*}
    \xymatrix{G_\alpha\times_{G_{\alpha-1}}X_{\alpha-1}
        \ar[d]^{\id\times g_{\alpha-1}} \ar[r]^-{i_\alpha} & \cyl(f_\alpha) \ar[d] \\
        G_\alpha\times_{G_{\alpha-1}}Y_{\alpha-1} \ar[r] & Y_\alpha
    }
\end{equation*}
in which $i_\alpha$ is the obvious inclusion into the mapping cylinder of
$f_\alpha$ and $g_{\alpha-1}$ the up to $G_{\alpha-1}$-homotopy unique homotopy
equivalence which comes from the universal property of
$\EGF{G_{\alpha-1}}{\calsfg}$. Hence, $Y_{\alpha}$ is $G_{\alpha}$-homotopy
equivalent to $X_{\alpha}$ and therefore a model for
$\EGF{G_{\alpha}}{\calsfg}$. Moreover, $Y_{\alpha}$ is clearly
$(n+1)$-dimensional. Finally, if $\alpha$ is a limit ordinal, we define
$Y_{\alpha}$ to be the union of the $G_{\alpha}\times_{G_{\beta}}Y_{\beta}$ for
$\beta<\alpha$.
\end{proof}

\begin{example} \label{exa:cor_hdimG(eub(G)_G_locally_finite}
Let $G$ be an infinite locally finite group of cardinality $\aleph_n$. Then
it is consistent with ZFC (Zermelo-Fraenkel plus axiom of choice) that 
\begin{equation*}
    \hdim^G(\eub{G}) = n+1.
\end{equation*}
The existence of an $(n+1)$-dimensional model for $\eub{G}$ follows from
Theorem~\ref{thm:locally_F_model} and has already been proved by
Dicks-Kropholler-Leary~\cite[Theorem~2.6]{Dicks-Kropholler-Leary-Thomas(2002)}.
Now the claim follows from the fact that it is consistent with ZFC that the cohomological dimension over $\IQ$
of $G$ satisfies $\cd_{\IQ}(G)=n+1$
(see~\cite[Theorem~A]{Kropholler-Thomas(1997)}).

We emphasize that the assumption in Example~\ref{exa:calf_is_fin_gen_sub} that
$G$ is countable is necessary. Namely,
Dunwoody~\cite[Theorem~1.1]{Dunwoody(1979)} has proved that for a group $K$
there is a one-dimensional model for $\eub{K}$ if and only if $\cd_{\IQ}(K)\le
1$.
\end{example}


\subsection{Low dimensions}
\label{subsec:low_dimensions}

For countable groups $G$, the question whether there are models for $\eub{G}$
and $\edub{G}$ of small dimension has the following answer. The assumption that
$G$ is countable is essential (see
Example~\ref{exa:cor_hdimG(eub(G)_G_locally_finite}).

\begin{theorem} \label{the:low_hdim}

\begin{enumerate}

\item \label{the:low_hdim:G_locally_virtually_cyclic}
Let $G$ be a countable group which is locally virtually cyclic. Then
\begin{equation*}
    \hdim^G(\eub{G}) =
    \begin{cases}
        0 & \text{if $G$ is finite;} \\
        1 & \text{if $G$ is infinite and either locally finite} \\
        & \hfill \text{or virtually cyclic;} \\
        2 & \text{otherwise,}
    \end{cases}
\end{equation*}
and
\begin{equation*}
    \hdim^G(\edub{G}) =
    \begin{cases}
        0 & \text{if $G$ is virtually cyclic;} \\
        1 & \text{otherwise;}
    \end{cases}
\end{equation*}

\item \label{the:low_hdim(eub(G)_le_1}
Let $G$ be a countable group satisfying $\hdim^G(\eub{G}) \le 1$. Then
\begin{equation*}
    \hdim^G(\edub{G}) =
    \begin{cases}
        0 & \text{if $G$ is virtually cyclic;} \\
        1 & \text{if $G$ is locally virtually cyclic but} \\
        & \hfill \text{not virtually cyclic;} \\
        2 & \text{otherwise.}
    \end{cases}
\end{equation*}
\end{enumerate}
\end{theorem}
\begin{proof}%
\ref{the:low_hdim:G_locally_virtually_cyclic} The claim about
$\hdim^G(\edub{G})$ follows from Example~\ref{exa:calf_is_fin_gen_sub}.

We conclude $\hdim^G(\eub{G}) \le 2$ from $\hdim^G(\edub{G}) \le 1$ and
Corollary~\ref{cor:dim(eub(G)_le_dim(edub(G)_plus_1}~\ref{cor:dim(eub(G)_le_dim(edub(G)_plus_1:mindim}.
If $G$ is locally finite, then $\hdim^G(\eub{G})\le 1$ has already been proved
in Example~\ref{exa:calf_is_fin_gen_sub}.  Let $G$ be a countable locally
virtually cyclic group and $\hdim^G(\eub{G})\le 1$. It remains to show that $G$
is virtually cyclic or locally finite. Assume that $G$ is not locally finite,
then we can choose choose a sequence of infinite virtually cyclic subgroups
$V_0\subseteq V_1 \subseteq V_2 \subseteq \ldots $ such that $G = \bigcup_{i
\ge 0} V_i$, and we must show that $G$ itself is virtually cyclic.

We begin with the case that each $V_i$ is an infinite virtually cyclic group of
type I.  Type I means that we can find an epimorphism $p_i \colon V_i \to C_i$
onto an infinite cyclic group $C_i$.  Let $F_i$ be the kernel of $p_i$, which
is a finite group. Thus we obtain a nested sequence of inclusions of finite
groups $F_1 \subseteq F_2 \subseteq \ldots$ and a sequence of inclusions of
infinite cyclic groups $C_1 \subseteq C_2 \subseteq \ldots$ such that there
exists a short exact sequence $1 \to F_i \to V_i \to C_i \to 1$ for each $i$.

Now we deal with the hardest step in the proof, where we show that the sequence
$F_0 \subseteq F_1 \subseteq F_2 \subseteq \ldots$ is eventually stationary.
Let $t \in V_0$ be an element which is mapped to a generator under the
epimorphism $p_0\colon V_0 \to C_0$. If we consider $t$ as an element in $V_i$,
conjugation with $t$ induces an automorphism $\phi_i\colon F_i \to F_i$.
Obviously, the restriction of $\phi_i$ to $F_{i-1}$ is $\phi_{i-1}$. We obtain
a nested sequence of virtually cyclic groups $F_0 \rtimes \IZ \subseteq F_1
\rtimes \IZ \subseteq F_2 \rtimes \IZ \subseteq \ldots$, where the $i$-th
semi-direct product $F_i\rtimes \IZ$ is to be understood with respect to
$\phi_i$. Since $K := \bigcup_{i \ge 0} F_i \rtimes \IZ$ is a subgroup of $G$,
it satisfies $\hdim^K(\eub{K}) \le 1$. This implies $H^2_K(\eub{K};M) = 0$ for
any $\IQ[K]$-module $M$.

The $\IQ[K]$-chain map $C_*(EK) \to C_*(\eub{K})$ of the cellular
$\IQ[K]$-chain complexes coming from a $K$-map $EK \to \eub{K}$ induces an
isomorphism on homology. Since both of these $\IQ[K]$-chain complexes are
projective, this map is a $\IQ[K]$-chain homotopy equivalence.  Hence, we
obtain for every $\IQ[K]$-module $M$ an isomorphism
\begin{multline*}
    H^k(K;M) = H^k(\hom_{\IQ[K]}(C_*(EK),M)) \cong
    H^k(\hom_{\IQ[K]}(C_*(\eub{K}),M)) \\
    = H^k_K(\eub{K};M).
\end{multline*}
In particular, $H^2(K;M) = 0$.

If we put $F = \bigcup_{i \ge 0} F_i$, then the collection of the $\phi_i$-s
yields an automorphism $\phi \colon F \to F$, with respect to which $K = F
\rtimes \IZ$. Let $M$ be any $\IQ[K]$-module.  The Hochschild-Serre spectral
sequence associated to $1 \to F \to K = F \rtimes \IZ \to \IZ \to 1$ yields the
exact sequence
\begin{equation*}
    H^1(F;M) \xrightarrow{\id - \phi_*} H^1(F;M) \to H^2(K;M) \to H^2(F;M)
\end{equation*}
where the automorphism $\phi_*$ of $H^1(F;M)$ is the one induced by $\phi
\colon F \to F$ and the $\IQ[F\rtimes \IZ]$-structure on $M$. Hence
\begin{equation*}
    H^2(K;M) \cong \coker\left(\id - \phi_* \colon H^1(F;M) \to H^1(F;M)\right)
\end{equation*}
since we will show in a moment that $H^2(F;M)=0$.

Next, we compute $H^*(F;M)$ for a given $\IQ[F]$-module $M$. We have already
explained the isomorphism
\begin{equation*}
    H^1(F;M) = H^1(\hom_{\IQ[F]}(C_*(EF),M)) \cong
    H^1(\hom_{\IQ[F]}(C_*(\eub{F}),M)).
\end{equation*}
It follows from Lemma~\ref{lemma:hocolim_G/im(psi_i)_simeq_pt} that we obtain
an $F$-pushout
\begin{equation*}
    \xymatrix@C=4em{\coprod_{i\ge 0} F/F_i \times S^0 \ar[d] \ar[r]^-{\coprod_{i\ge 0} k_i}
        & \coprod_{i\ge 0} F/F_i \ar[d] \\
        \coprod_{i\ge 0} F/F_i \times D^1 \ar[r] & \eub{F}
    }
\end{equation*}
where the left vertical arrow is the canonical inclusion and $k_i \colon F/F_i
\times S^0 \to \coprod_{i\ge 0} F/F_i$ sends $(u,-1)$ to $u \in F/F_i$ and
$(u,1)$ to the image of $u$ under the projection $F/F_i \to F/F_{i+1}$. Hence,
the cellular $\IQ[F]$-chain complex of $\eub{F}$ is concentrated in dimensions
$0$ and $1$ and is given there by
\begin{equation*}
    \bigoplus_{i\ge 0} \IQ[F/F_i] \xrightarrow{\bigoplus_{i\ge 0} l_i}
    \bigoplus_{i\ge 0} \IQ[F/F_i],
\end{equation*}
where $l_i$ sends $u \in \IQ[F/F_i]$ to the difference of the element given by
$u$ itself and the element given by the image of $u$ under $\IQ[\pr_i] \colon
\IQ[F/F_i] \to \IQ[F/F_{i+1}]$ for $\pr_i \colon F/F_i \to F/F_{i+1}$ the
canonical projection. Therefore, $H^k(F;M)= 0$ for $k \ge 2$, whereas
$H^1(F;M)$ is the cokernel of the $\IQ[F]$-homomorphism
\begin{equation*}
    \delta_M \colon \prod_{i \ge 0} M^{F_i} \to \prod_{i \ge 0} M^{F_i}, \quad
    (x_i)_{i \ge 0} \mapsto (x_i - x_{i+1})_{i \ge 0}.
\end{equation*}
Let us show that $H^1(F;\IQ[F])= 0$ if and only if the sequence $F_0 \subseteq
F_1 \subseteq F_2 \subseteq \ldots$ is eventually stationary.  If it becomes
stationary, then $F$ is finite, so $H^1(F;\IQ[F])= 0$. Suppose that
$H^1(F;\IQ[F])=0$ and let $N_{F_i} \in \IQ[F]^{F_i}$ be the element $\sum_{f
\in F_i} f$.  Since $\delta_{\IQ[F]}$ is surjective, there exists $(x_i)_{i \ge
0} \in \prod_{i \ge 0} \IQ[F]^{F_i}$ satisfying $N_{F_i} = x_i -x_{i+1}$. Thus,
we get for $i \ge 1$
\begin{equation*}
    x_i = x_0 - N_{G_{i-1}} - N_{G_{i-2}} - \cdots - N_{G_0}.
\end{equation*}
After choosing $i_0$ in such a way that $x_0 \in \IQ[F_{i_0}]$, we get that
$x_i \in \IQ[F_{i-1}]$ holds for for $i > i_0$. Since $x_i \in \IQ[F]^{F_i}$,
we conclude that $F_i = F_{i-1}$ holds for $i > i_0$, which means that $F_0
\subseteq F_1 \subseteq F_2 \subseteq \ldots$ is eventually stationary. For
more information about $H^*(F;\IQ[F])$ for locally finite groups $F$, we refer
for instance to~\cite{Dicks-Kropholler-Leary-Thomas(2002)}.

Now suppose that $M$ is a $\IQ[K]$-module. We obtain the commutative diagram
\begin{equation*}
    \xymatrix{\prod_{i \ge 0} M^{F_i} \ar[d]_{\id - l_t} \ar[r]_{\delta_M}
        & \prod_{i \ge 0} M^{F_i} \ar[d]_{\id - l_t} \ar[r]
        & H^1(F;M) \ar[d]_{\id - \phi_*} \ar[r] & 0 \\
        \prod_{i \ge 0} M^{F_i} \ar[r]_{\delta_M}
        & \prod_{i \ge 0} M^{F_i} \ar[r] & H^1(F;M) \ar[r] & 0
    }
\end{equation*}
where $l_t \colon M \to M$ is multiplication with the generator of $\IZ$
considered as an element in $K= F \rtimes \IZ$.

Specializing to $M = \IQ[K] = \IQ[F \rtimes \IZ]$, we can write down the
following exact sequence
\begin{equation*}
    0 \to \IQ[F\rtimes \IZ] \xrightarrow{\id - l_t} \IQ[F\rtimes \IZ]
    \xrightarrow{\epsilon} \IQ[F] \to 0,
\end{equation*}
where $\epsilon$ sends $\sum_{n \in \IZ} r_i \cdot f_it^i$ to $\sum_{n \in \IZ}
r_i \cdot \phi^{-i}(f_i)$ for $r_i \in \IQ$, $f_i \in F$ and $t \in \IZ$ a
fixed generator. The composition
\begin{equation*}
    \prod_{i \ge 0} \IQ[F \rtimes \IZ]^{F_i} \xrightarrow{\delta_{\IQ[F \rtimes \IZ]}}
    \prod_{i \ge 0} \IQ[F \rtimes \IZ]^{F_i} \xrightarrow{\epsilon} \prod_{i \ge 0}
    \IQ[F]^{F_i}
\end{equation*}
agrees with the composition
\begin{equation*}
    \prod_{i \ge 0} \IQ[F \rtimes \IZ]^{F_i} \xrightarrow{\epsilon} \prod_{i \ge 0}
    \IQ[F]^{F_i} \xrightarrow{\delta_{\IQ[F]}} \prod_{i \ge 0} \IQ[F]^{F_i}
\end{equation*}
since $\epsilon\circ l_t = \epsilon$. Hence we obtain an isomorphism
\begin{multline*}
    H^1(F;\IQ[F]) \cong \coker\bigl(\delta_{\IQ[F]}\colon H^1(F;\IQ[F]) \to
    H^1(F;\IQ[F])\bigr) \\
    \cong \coker\bigl(\id - \phi_* \colon H^1(F;\IQ[F \rtimes \IZ]) \to
    H^1(F;\IQ[F \rtimes \IZ])\bigr) \cong H^2(K;\IQ[K]).
\end{multline*}
As $H^2(K;\IQ[K])$ is trivial, the same holds for $H^1(F;\IQ[F])$. We have
already shown that this implies that $F = \bigcup_{i \ge 0} F_i$ is finite.

Returning to the beginning of the proof and denoting by $A$ the colimit of the
system $C_0 \subseteq C_1 \subseteq C_2 \subseteq \ldots$ of infinite cyclic
groups, we obtain the exact sequence $1\to F\to G\to A\to 1$. Since
$\hdim^G(\eub{G}) \le 1$ by assumption, the cohomological dimension of $G$ over
$\IQ$ satisfies $\cd_{\IQ}(G)\le 1$. Consequently, also $\cd_{\IQ}(A)\le 1$ as
$F$ is finite. By a result of Dunwoody~\cite[Theorem~1.1]{Dunwoody(1979)}
already mentioned in Example~\ref{exa:cor_hdimG(eub(G)_G_locally_finite}, this
means that $\hdim^A(\eub{A}) \le 1$. However, as $A$ is torsionfree, this
implies that $A$ is free (see~\cite{Stallings(1968)} and~\cite{Swan(1969)}).
Since $\IZ$ is the only non-trivial abelian free group, $A$ is infinite cyclic,
so $G$ is virtually cyclic.  This finishes the proof in the case that each
$V_i$ is infinite cyclic of type I.

We are left to treat the case that there exists $i_0$ such that $V_{i_0}$ is
infinite virtually cyclic of type {II}, i.e., there exists an epimorphism
$V_{i_0}\to D_\infty$ onto the infinite dihedral group. Since $V_{i_0}$ is a
subgroup of $V_i$ for $i \ge i_0$, the virtually cyclic group $V_i$ is of type
{II} for $i \ge i_0$ or, equivalently, $H_1(V_i)$ is finite for $i\ge i_0$. Let
$W_i$ be the commutator of $V_i$. Then we obtain a sequence of virtually cyclic
subgroups $W_0\subseteq W_1 \subseteq W_2 \subseteq \ldots$ such that $W_i$ is
of type I and has finite index in $V_i$ for $i \ge i_0$. Let $W$ be the
subgroup $\bigcup_{i \ge 0} W_i$ of $G$.  Since $\dim^G(\eub{G}) \le 1$ by
assumption, we have $\hdim^W(\eub{W}) \le 1$ as well. Thus, by what we have
already proved above, $W$ is virtually cyclic and, in particular, finitely
generated. Hence, the sequence $W_0 \subseteq W_1 \subseteq W_2 \subseteq
\ldots$ is eventually stationary, and we can assume without loss of generality
that $W = W_i$ for all $i \ge 0$. This implies that each map $H_1(V_i) \to
H_1(V_{i+1})$ is injective. Putting $A := \bigcup_{i \ge 0} H_1(V_i)$, we
obtain an exact sequence $1 \to W \to G \to A \to 0$, where $W$ is an infinite
virtually cyclic group of type I. We are going to show that $A$ is finite.

Let $K$ denote the kernel of the canonical epimorphism $W \to C$ for $C =
H_1(W)/\tors(H_1(W))$. Then $K$ is a finite characteristic subgroup of $W$ and
hence normal in $G$, so there is an exact sequence $1 \to K \to G \to G/K \to
1$. Note that we have already explained above (using
\cite[Theorem~1.1]{Dunwoody(1979)}) why $\hdim^G(\eub{G}) \le 1$ implies
$\hdim^{G/K}(\eub{(G/K)}) \le 1$. Furthermore, there is an exact sequence $1\to
C\to G/K\to A\to 1$. Let $A \to \aut(C)$ be the associated conjugation
homomorphism, let $A' \subseteq A$ be its kernel, and let $G'$ the preimage of
$A'$ under the epimorphism $G/K \to A$. Then we obtain a central extension $1
\to C \to G' \to A' \to 1$ such that $\hdim^{G'}(\eub{G'}) \le 1$ and $A'
\subseteq A$ is a subgroup of finite index.  We have $A' = \bigcup_{i \ge 0}
(A' \cap H_1(V_i))$. The computations appearing in the proof above imply that
either $A'$ is a finite abelian group or $H^1(A';\IQ[A']) \cong
\coker(\delta_{\IQ[A']})\not= 0$, and that $H^k(A';\IQ[A']) = 0$ for $k \ge 2$.
From the Hochschild-Serre cohomology spectral sequence applied to $1 \to C \to
G' \to A' \to 1$ we obtain an isomorphism $H^1(A';\IQ[A']) \cong
H^2(G';\IQ[A'])$, and the latter term vanishes as $\hdim^{G/K}(\eub{(G/K)}) \le
1$. We conclude that $A'$ and hence $A$ are finite, showing that $G$ is
virtually cyclic.  This finishes the proof of
assertion~\ref{the:low_hdim:G_locally_virtually_cyclic}.
\\[1mm]\ref{the:low_hdim(eub(G)_le_1}
First we show that $\hdim^G(\edub{G}) \le 2$ provided that $\hdim^G(\eub{G})
\le 1$ and $G$ is countable.  Consider an infinite virtually cyclic subgroup
$V\subseteq G$.  Let $H\subseteq N_G[V]$ be a finitely generated subgroup such
that $V \subseteq H$.  From $\hdim^G(\eub{G}) \le 1$ we get $\hdim^H(\eub{H})
\le 1$, and since $H$ is finitely generated, this implies that it is virtually
finitely generated free (see \cite[Theorem
1]{Karrass-Pietrowski-Solitar(1973)}). In particular, $H$ is word-hyperbolic,
so we conclude from Theorem~\ref{the:max_vc_sg} and Example~\ref{exa:word-hyp}
that $N_H[V]$ is virtually cyclic.  Since $H = H \cap N_G[V] = N_H[V]$, this
means that $H$ is virtually cyclic. Thus, we have shown that $N_G[V]$ is
locally virtually cyclic, and as $\hdim^G(\eub{G}) \le 1$ implies
$\hdim^{N_G[V]}(\eub{N_G[V]}) \le 1$,
assertion~\ref{the:low_hdim:G_locally_virtually_cyclic} shows that $N_G[V]$ is
locally finite or virtually cyclic. In fact, $N_G[V]$ contains the infinite
virtually cyclic group $V$, so $N_G[V]$ is virtually cyclic. This implies that
$\calvcyc[V]$ consists of all subgroups of $N_G[V]$ and, hence, $\pt$ is a
model for $\EGF{N_G[V]}{\calvcyc[V]}$. Summarizing, we can arrange that
$\dim(\eub{N_G[V]}) \le 1$ and $\dim(\EGF{N_G[V]}{\calvcyc[V]}) \le 1$ for
every infinite virtually cyclic subgroup $V$ of $G$. Now
Theorem~\ref{the:passing_from_calf_to_calg} together with
Remark~\ref{rem:how_to_apply_theorem_the:passing_from_calf_to_calg} implies
$\hdim^G(\edub{G}) \le 2$.

If $G$ is countable and locally virtually cyclic, then $\hdim^G(\edub{G}) \le 1
$ has already been proved in
assertion~\ref{the:low_hdim:G_locally_virtually_cyclic}. It remains to show
that $G$ is locally virtually cyclic provided that $\hdim^G(\eub{G})$ and
$\hdim^G(\edub{G})$ are less or equal to $1$.  Let $G_0 \subseteq G$ be a
finitely generated subgroup of $G$.  Then $\hdim^G(\eub{G_0})$ and
$\hdim^G(\edub{G_0})$ are less or equal to $1$.  This implies that $G_0$
contains a finitely generated free subgroup $G_1$ of finite index (see
\cite[Theorem 1]{Karrass-Pietrowski-Solitar(1973)}). The group $G_1$ satisfies
$N\!M_{\calfin\subseteq\calvcyc}$ by Example~\ref{exa:word-hyp}, so from
Corollary~\ref{cor:Efin_and_Evcyc_model_with_N_GM_is_M} we get a cellular
$G_1$-pushout
\begin{equation*}
    \xymatrix{\coprod_{V\in\calm} G_1\times_{V}\eub{V}
        \ar[d]^{\coprod_{V\in\calm} \id_{G_1}\times_{V} i_V} \ar[r]^-i
        & \eub{G_1} \ar[d] \\
        \coprod_{V\in\calm} G_1\times_{V} \cone(\eub{V}) \ar[r] & \edub{G_1}
    }
\end{equation*}
where $\calm$ is a complete system of representatives of the conjugacy classes
of maximal infinite virtually cyclic subgroups of $G_1$ and $i_V$ is the
obvious inclusion. Dividing out the $G_1$-action and then taking the associated
Mayer-Vietoris exact sequence yields an injection
\begin{equation*}
    \bigoplus_{V \in \calm} H_1(V\backslash\eub{V}) \to H_1(G_1\backslash\eub{G_1})
\end{equation*}
because $\hdim^{G_1}(\edub{G_1})\le 1$. Since $G_1$ is finitely generated free,
$H_1(G_1\backslash\eub{G_1})$ is finitely generated free and
$H_1(V\backslash\eub{V})$ is isomorphic to $\IZ$ for every $V \in \calm$. This
implies that $\calm$ is finite, i.e., $G_1$ contains only finitely many
infinite cyclic subgroups up to conjugacy. Thus, $G_1$ is either trivial or
infinite cyclic, showing that $G_0$ is virtually cyclic. Hence $G$ is locally
virtually cyclic. This finishes the proof of Theorem~\ref{the:low_hdim}.
\end{proof}

\begin{example} \label{exa:hdim_for_eub_and_edub_of_Z[1/p]}
For the prime number $p$, let $\IZ[1/p]$ be the subgroup of $\IQ$ consisting of
rational numbers $x \in \IQ$ for which $p^n \cdot x \in \IZ$ for some positive
integer $n$ (see Example~\ref{exa:N_G[V]_is_N_GC}). Every finitely generated
subgroup of $\IZ[1/p]$ is trivial or infinite cyclic. Thus,
Theorem~\ref{the:low_hdim} implies
\begin{align*}
    \hdim^{\IZ[1/p]}(\eub{\IZ[1/p]}) = 2; \\
    \hdim^{\IZ[1/p]}(\edub{\IZ[1/p]}) = 1. \\
\end{align*}
\end{example}

\begin{remark} \label{rem:characterization_of_countable_locally_finitw_groups}
Let $G$ be a countable group. Then Theorem~\ref{the:low_hdim} implies that $G$
is infinite locally finite if and  only if $\hdim^G(\eub{G}) =
\hdim^G(\edub{G}) = 1$ holds.
\end{remark}


\typeout{-------------------- Section 6 --------------------------}

\section{Equivariant homology of relative assembly maps}
\label{sec:Equivariant_homology_of_relative_assembly_maps}

Let $\calh^?_*$ be an equivariant homology theory in the sense
of~\cite[Section~1]{Lueck(2002b)}. Our main example is the equivariant homology
theory $H_*^?(-;\bfK_R)$ appearing in the $K$-theoretic Farrell-Jones
conjecture, where $R$ is a ring (here, rings are always assumed to be
associative and unital). It has the property that
$$H_n^G(G/H;\bfK_R) = H_n^H(\pt;\bfK_R) = K_n(RH)$$
holds for every subgroup $H \subseteq G$, where $K_n(RH)$ denotes the algebraic
$K$-theory of the group ring $RH$. The Farrell-Jones conjecture for a group $G$
and a ring $R$ says that the assembly map, which is the map induced by the
projection $\edub{G} \to G/G$, is an isomorphism
\begin{equation*}
    H_n^G(\edub{G};\bfK_R) \xrightarrow{\cong} H_n^G(G/G;\bfK_R) = K_n(RG)
\end{equation*}
for all $n \in\IZ$.

Bartels~\cite{Bartels(2003b)} has shown that for every group $G$, every ring $R$, and
every $n \in \IZ$ the relative assembly map
\begin{equation*}
    H_n^G(\eub{G};\bfK_R) \to H_n^G(\edub{G};\bfK_R)
\end{equation*}
is split-injective. Hence, the source of the assembly map appearing in the
Farrell-Jones conjecture can be computed in two steps, involving the
computation of $H_n^G(\eub{G};\bfK_R)$ and the computation of the remaining
term which we denote by $H_n^G(\eub{G}\to\edub{G};\bfK_R)$. Rationally,
$H_n^G(\eub{G};\bfK_R)$ can be computed using equivariant Chern characters
(see~\cite[Section~1]{Lueck(2002b)}) or with the help of nice models for
$\eub{G}$.

The goal of this  section is to give some information about $H_n^G (
\eub{G} \to \edub{G}; \bfK_R )$. Namely, using the induction structure on
$\calh^?_*$, Corollary~\ref{cor:Efin_and_Evcyc_model_assuming_(McalG,calF)}
implies:

\begin{corollary} \label{cor:calh_assuming_(McalG,calF)}
Let $\calh^?_*$ be an equivariant homology theory.  Suppose that the group $G$
satisfies $(M_{\calfin \subseteq \calvcyc})$. Let $\calm$ be a complete system
of representatives of the conjugacy classes of maximal infinite virtually
cyclic subgroups of $G$. Then, for $n \in \IZ$, we obtain an isomorphism
\begin{equation*}
    \bigoplus_{V \in \calm} \calh_n^{N_GV}(\eub{N_GV} \to EW_GV)
    \xrightarrow{\cong} \calh_n^G ( \eub{G} \to \edub{G}).
\end{equation*}
\end{corollary}

\begin{remark}\label{rem:Term_calh_nN_GVeub(N_GV)_to_EW_GV)}
The term $\calh_n^{N_GV}(\eub{N_GV} \to EW_GV)$ appearing in
Corollary~\ref{cor:calh_assuming_(McalG,calF)} can be analyzed further.

Namely, by assigning to a free $W_GV$-$CW$-complex $X$ the $\IZ$-graded abelian
group $\calh^{N_GV}_*(\eub{N_GV} \times X \to X)$, we obtain a $W_GV$-homology
theory for free $W_GV$-$CW$-complexes. Here, we consider $X$ as an
$N_GV$-$CW$-complex by restriction with the projection $N_GV \to W_GV$, equip
$\eub{N_GV}\times X$ with the diagonal $N_GV$-action and let $\eub{N_GV}\times
X \to X$ be the projection onto $X$. There is an Atiyah-Hirzebruch spectral
sequence converging to $\calh^{N_GV}_{p+q}(\eub{N_GV} \times X \to X)$ whose
$E^2$-term is
\begin{equation*}
    E^2_{p,q} = H^{W_GV}_p(EW_GV;\calh_q^{N_GV}(EN_GV \times W_GV \to W_GV)),
\end{equation*}
where the right $\IZ[W_GV]$-module structure on $\calh_q^{N_GV}(\eub{N_GV}
\times W_GV \to W_GV)$ comes from the obvious right $W_GV$-action on $W_GV$ and
the trivial $W_GV$-action on $\eub{N_GV}$.

For any $N_GV$-space $Y$, the map $N_GV \times_V \res_{N_GV}^VY \to Y \times
W_GV$ given by sending $(g,y)$ to $(gy,gV)$ is an $N_GV$-homeomorphism. In
particular, since $\res_{N_GV}^V\eub{N_GV} $ is a model for $\eub{V}$, we
obtain an $N_GV$-homeomorphism $N_GV \times_V \eub{V} \to \eub{N_GV} \times
W_GV$.  Together with the induction structure, this induces an isomorphism of
abelian groups
\begin{equation*}
    \calh_q^{N_GV}(EN_GV \times W_GV \to W_GV)) \cong \calh_q^{V}(\eub{V} \to \pt),
\end{equation*}
which becomes an isomorphism of $\IZ[N_GV]$-modules if we equip the target with
the following $W_GV$-action.

Given $\overline{g} \in W_GV$, choose a preimage $g \in N_GV$ under the
projection $N_GV \to W_GV$.  Conjugation with $g$ yields a group homomorphism
$c(g) \colon V \to V$. The induction structure yields an isomorphism
$\calh^V_q(\eub{V} \to \pt) \to \calh^V_q(\ind_{c(g)}\eub{V} \to
\ind_{c(g)}\pt)$. Let $f \colon \ind_{c(g)}\eub{V} \to \eub{V}$ be the
$V$-homotopy equivalence which is unique up to $V$-homotopy. Obviously,
$\ind_{c(g)}\pt = \pt$. Hence, $f$ induces a homomorphism $f_q \colon
\calh_q^V(\ind_{c(g)}\eub{V} \to \ind_{c(g)}\pt) \to \calh^V_q(\eub{V} \to
\pt)$.  The composition of these two homomorphisms is an automorphism of
$\calh^V_q(\eub{V} \to \pt)$, which we define to be multiplication with
$\overline{g}$. It is easy to check that this definition is independent of the
choice of the preimage $g$ of $\overline{g}$, and that it defines a
$W_GV$-action on $\calh_q^{V}(\eub{V} \to \pt)$ which is compatible with the
one on $\calh_q^{N_GV}(EN_GV \times W_GV \to W_GV)$.
\end{remark}

\begin{example} \label{exa:rem:Term_calh_nN_GVeub(N_GV)_to_EW_GV;bkk_R)}
Let us consider the special case that the equivariant homology theory is
$H^?_*(-;\bfK_R)$, the one appearing in the $K$-theoretic Farrell-Jones
conjecture. Let $V$ be an infinite virtually cyclic group. If $V$ is of type
$I$, we can write $V$ as a semi-direct product $F \rtimes \IZ$, and
$H^V_*(\eub{V} \to \pt;\bfK_R)$ is the non-connective version of Waldhausen's
Nil-term associated to this semi-direct product
(see~\cite[Section~9]{Bartels-Lueck(2006)}). If $V$ is of type {II}, then it
can be written as an amalgamated product $V_1 \ast_{V_0} V_2$ of finite groups,
where $V_0$ has index two in both $V_1$ and $V_2$. In this case, $H^V_*(\eub{V}
\to \pt;\bfK_R)$ is the non-connective version of Waldhausen's
Nil-term associated to this amalgamated product
(see~\cite[Section~9]{Bartels-Lueck(2006)}). The $W_GV$-action on these
Nil-terms comes from the action of $N_GV$ on $V$ by conjugation and the fact
that inner automorphisms of groups induce the identity on algebraic $K$-groups
associated to group rings.

In the case of $L$-theory, the terms $H^V_*(\eub{V} \to
\pt;\bfL_R^{\langle-\infty\rangle})$ vanish if $V$ is of type I and are given
by UNil-terms if $V$ is of type {II} (see~\cite[Lemma~4.2]{Lueck(2005heis)}).
\end{example}

\begin{example} \label{exa:calh_assuming_(NMcalG,calF)}
Let $\calh^?_*$ be an equivariant homology theory.  Suppose that the group $G$
satisfies $(N\!M_{\calfin \subseteq \calvcyc})$. Let $\calm$ be a complete system
of representatives of the conjugacy classes of maximal infinite virtually
cyclic subgroups of $G$. Then, for $n\in\IZ$,
Corollary~\ref{cor:calh_assuming_(McalG,calF)} yields an isomorphism
\begin{equation*}
    \bigoplus_{V \in \calm} \calh_n^{V}(\eub{V} \to \pt)
    \xrightarrow{\cong} \calh_n^G(\eub{G} \to \edub{G}).
\end{equation*}
Now assume in addition that $G$ is torsionfree. If the equivariant homology
theory is $H^?_*(-,\bfK_R)$, then $H^V_*(\eub{V} \to \pt;\bfK_R)$ reduces to
$N\!K_n(R) \oplus N\!K_n(R)$. Here, $N\!K_n(R)$ is the $n$-th Bass-Nil-group
which is defined as the cokernel of the obvious split injection $K_n(R) \to
K_n(R[t])$. In the case of $L$-theory, we get $H^V_*(\eub{V} \to
\pt;\bfL_R^{\langle-\infty\rangle}) = 0$.
\end{example}


\typeout{-------------------- References -------------------------------}

\def\cprime{$'$} \def\polhk#1{\setbox0=\hbox{#1}{\ooalign{\hidewidth
  \lower1.5ex\hbox{`}\hidewidth\crcr\unhbox0}}}


\end{document}